\documentclass[11pt,twoside,a4paper,leqno]{amsart}
\calclayout
\usepackage{amssymb}
\newif\ifaddpics\addpicstrue% Comment out \addpicstrue to remove pictures.
\ifaddpics\usepackage{graphicx}\fi
\usepackage{diagrams}
\makeatletter
\def\swappedhead@plain#1#2#3{%
  \thmnumber{\@upn{\mdseries #2}}\thmname{\@ifnotempty{#2}{. }#1}%
  \thmnote{ \textmd{\upshape(#3)}}}
\makeatother
\theoremstyle{plain}
\newtheorem{result}{Theorem}

\swapnumbers
\newtheorem{thm}[subsubsection]{Theorem}
\newtheorem{prop}[subsubsection]{Proposition}
\newcommand{\tcite}[1]{\textup{\cite{#1}}}      % Theorem references
\numberwithin{equation}{section}

\newcommand{\acknowledge}{\subsection*{Acknowledgements}}
\newcommand{\thismonth}{\ifcase\month\or
  January\or February\or March\or April\or May\or June\or
  July\or August\or September\or October\or November\or December\fi
  \space\number\year}
\makeatletter
\newcommand{\low}{\@ifnextchar^{}{^{\vphantom x}}}
\newcommand{\high}{\@ifnextchar_{}{_{\vphantom I}}}
\makeatother
\DeclareSymbolFont{script}{U}{eus}{m}{n}
\DeclareSymbolFontAlphabet{\mathscr}{script}
\DeclareMathSymbol{\EuWedge}{0}{script}{"5E}
\DeclareMathAlphabet{\mathrmsl}{OT1}{cmr}{m}{sl}
\newcommand{\rssymb}[2]{\newcommand{#1}{{\mathrmsl{#2}}}}
\newcommand{\calsymb}[2]{\newcommand{#1}{{\mathcal{#2}}}}
\newcommand{\bbsymb}[2]{\newcommand{#1}{{\mathbb{#2}}}}
\newcommand{\lieoper}[2]{\newcommand{#1}{\mathop
  {\mathfrak{#2}\null}\nolimits}}
\newcommand{\oper}[3][n]{\newcommand{#2}{\mathop
  {\mathrm{#3}\null}\ifx n#1\nolimits\else\limits\fi}}
\newcommand{\rsoper}[3][n]{\newcommand{#2}{\mathop
  {\mathrmsl{#3}\null}\ifx n#1\nolimits\else\limits\fi}}
\bbsymb\C{C} \bbsymb\F{F} \bbsymb\IH{H}\bbsymb\N{N} \bbsymb\IP{P}
\bbsymb\Q{Q} \bbsymb\R{R} \bbsymb\U{U} \bbsymb\V{V} \bbsymb\W{W} \bbsymb\Z{Z}
\calsymb\cA{A} \calsymb\cB{B} \calsymb\cC{C} \calsymb\cD{D} \calsymb\cE{E}
\calsymb\cF{F} \calsymb\cG{G} \calsymb\cH{H} \calsymb\cI{I} \calsymb\cJ{J}
\calsymb\cK{K} \calsymb\cL{L} \calsymb\cM{M} \calsymb\cN{N} \calsymb\cO{O}
\calsymb\cP{P} \calsymb\cQ{Q} \calsymb\cR{R} \calsymb\cS{S} \calsymb\cT{T}
\calsymb\cU{U} \calsymb\cV{V} \calsymb\cW{W} \calsymb\cX{X} \calsymb\cY{Y}
\calsymb\cZ{Z}
\newcommand{\eps}{\varepsilon}
\newcommand{\ve}{\varepsilon}
 
\newcommand{\lam}{\lambda}
\renewcommand{\geq}{\geqslant} \renewcommand{\leq}{\leqslant}
\oper\End{End} \oper\Hom{Hom}                    % Vector space constructions
\oper\Sym{Sym} \oper\Skew{Skew}
\oper\Aut{Aut}                                   % Group constructions
\oper\GL{GL} \oper\SL{SL}\oper\Symp{Sp}
\oper\CO{CO} \oper\On{O} \oper\SO{SO} \oper\Pin{Pin} \oper\Spin{Spin}
\oper\CU{CU} \oper\Un{U} \oper\SU{SU}
\rsoper\Diff{Diff} \rsoper\SDiff{SDiff}
\lieoper\der{der}                                % Lie algebra constructions
\lieoper\gl{gl} \lieoper\sgl{sl}\lieoper\symp{sp}
\lieoper\co{co} \lieoper\so{so} \lieoper\spin{spin}
\lieoper\cu{cu} \lieoper\un{u}  \lieoper\su{su}
\rsoper\Vect{Vect} \rsoper\Ham{Ham}
\newcommand{\ip}[1]{\langle#1\rangle}

\newcommand{\norm}[2][]{|\mkern-2mu|#2|\mkern-2mu|
  _{\lower1pt\hbox{${}_{#1}$}}}
\newcommand{\Norm}[2][]{\bigl|\mkern-3mu\bigr|#2\bigr|\mkern-3mu\bigr|
  _{\lower1pt\hbox{${}_{#1}$}}}

\newcommand{\restr}[1]{|_{#1}\low}

\newcommand{\cross}{\mathbin{{\times}\!}\low}

\newcommand{\mult}{^{\scriptscriptstyle\times}}

\newcommand{\tens}{\otimes}                 % small tensor product
\newcommand{\setdif}{\smallsetminus}%{\mathbin{\mathrm -}}% set difference
\newcommand{\connect}{\#}                   % connected sum
\newcommand{\del}{\partial}                 % directional derivative
\newcommand{\Proj}{\mathrmsl{P}}            % projective
\newcommand{\RP}[1]{\R\Proj^{#1}}           % real projective n-space
\newcommand{\CP}[1]{\C\Proj^{#1}}           % complex projective n-space
\newcommand{\Cinf}{\mathrm{C}^\infty}       % smooth sections
\newcommand{\ie}{\textit{i.e.}}             % \ie,
\rsoper\dimn{dim}                           % dimension
\rsoper\grad{grad}                          % gradient
\rsoper\kernel{ker}\rsoper\image{im}        % kernel and image
\rsoper\alt{alt}   \rsoper\sym{sym}         % alternating and symmetric part
\rsoper\Ad{Ad}     \rsoper\ad{ad}           % adjoint action or bundle
\rsoper\CoAd{CoAd} \rsoper\coad{coad}       % coadjoint action
\rsoper\trace{tr}  \rsoper\trfree{tf}       % trace and tracefree part
\rsoper\detm{det}                           % determinant
\rsoper\Vol{Vol}                            % volume
\rsoper\divg{div}                           % divergence
\rsoper\sign{sign}
\rssymb\iden{id}                            % identity
\rssymb\vol{vol}                            % volume element
\oper\Imag{Im}
\newcommand{\trho}{\tilde{\rho}}
\newcommand{\oH}{\overline{\cH}{}}
\newcommand{\oX}{\overline{X}}
\newcommand{\teta}{\tilde{\eta}}
\newcommand{\ty}{\tilde{y}}

\newcommand{\sd}{{\raise1pt\hbox{$\scriptscriptstyle +$}}}
\newcommand{\asd}{{\raise1pt\hbox{$\scriptscriptstyle -$}}}
\newcommand{\sdasd}{{\raise1pt\hbox{$\scriptscriptstyle\pm$}}}
\newcommand{\asdsd}{{\raise1pt\hbox{$\scriptscriptstyle\mp$}}}

\rsoper\scal{scal}
\newcommand{\Quabla}{\pmb{\square}}
\newbox\dotbox
\newcommand\dotify[1]{\setbox\dotbox=\hbox{${#1}$}%
  \setbox0=\hbox{${\kern -.5\wd\dotbox\relax #1}$}\wd0=0pt\box0}
\newcommand{\isot}[2]{\rLine^{\scriptstyle#1}_{\scriptstyle #2}}
\begin{document}
\title{Einstein metrics and complex singularities}
\author{David M. J. Calderbank}
\address{Department of Mathematics and Statistics\\
University of Edinburgh\\ King's Buildings, Mayfield Road\\ Edinburgh
EH9 3JZ\\ Scotland.}
\email{davidmjc@maths.ed.ac.uk}
\author{Michael A. Singer}
\email{michael@maths.ed.ac.uk}
\date{\thismonth}
\begin{abstract}
This paper is concerned with the construction of special metrics on
non-compact $4$-manifolds which arise as resolutions of complex orbifold
singularities. Our study is close in spirit to the construction of the
hyperk\"ahler gravitational instantons, but we focus on a different class
of singularities. We show that any resolution $X$ of an isolated cyclic
quotient singularity admits a complete scalar-flat K\"ahler metric (which
is hyperk\"ahler if and only if $c_1(X)=0$), and that if $c_1(X)<0$ then
$X$ also admits a complete (non-K\"ahler) self-dual Einstein metric of
negative scalar curvature. In particular, complete self-dual Einstein
metrics are constructed on simply-connected non-compact $4$-manifolds with
arbitrary second Betti number.

Deformations of these self-dual Einstein metrics are also constructed: they
come in families parameterized, roughly speaking, by free functions of one
real variable.

All the metrics constructed here are {\em toric} (that is, the isometry
group contains a $2$-torus) and are essentially explicit.  The key to the
construction is the remarkable fact that toric self-dual Einstein metrics
are given quite generally in terms of {\em linear} partial differential
equations on the hyperbolic plane.
\end{abstract}

\maketitle

\section{Introduction and main theorems}

If $\Gamma$ is a finite subgroup of $SU(2)$, then the complex orbifold
$\C^2/\Gamma$ has a canonical resolution $X$ with $c_1(X)=0$. This
non-compact complex surface carries a family of asymptotically locally
euclidean (ALE) hyperk\"ahler metrics, the so-called gravitational
instantons of Gibbons, Hawking, Hitchin and
Kronheimer~\cite{GiHa:gmi,Hit:pg,Kro:ale}.  In this paper, we extend this
picture by looking for `optimal' metrics on complex resolutions of other
surface singularities. Our methods apply to finite cyclic subgroups $\Gamma
\subset U(2)$ with the property that $\C^2/\Gamma$ has an isolated singular
point, the image of the origin in $\C^2$. These varieties are \emph{toric}:
there is a $\C\mult\cross\C\mult$ subgroup of $\GL_2(\C)$ commuting with
$\Gamma$, which acts on $\C^2/\Gamma$ and the resolution $X$.  If
$\Gamma\not\subset SU(2)$ then $c_1(X)\neq 0$ and so $X$ cannot carry a
hyperk\"ahler metric. However we shall find ALE K\"ahler metrics with zero
scalar curvature, and if $c_1(X)$ is negative definite, complete
asymptotically (locally) hyperbolic self-dual Einstein metrics (with
respect to the opposite orientation of $X$).

The situation is simplest if $\Gamma$ acts by scalar multiples of the
identity on $\C^2$. Then if $|\Gamma|=p$, $X$ is the total space of the
complex line bundle $\cO(-p)\to \CP{1}$ and for each $p$ (including $p=1$,
the blow up of $\C^2$) there is a $U(2)$-invariant ALE scalar-flat K\"ahler
(SFK) metric on $X$: this is due to Burns if $p=1$, to Eguchi--Hanson if
$p=2$, and to LeBrun for $p>2$ (see~\cite{LeBr:pac}). Note that among these
metrics, only Eguchi--Hanson is hyperk\"ahler, corresponding to $\Gamma =
\{\pm 1\}\subset SU(2)$.  However, the Burns metric is conformal to the
Fubini--Study metric on the punctured projective plane
$\CP2\setdif\{\infty\}$, which is a self-dual Einstein (SDE) metric of
positive scalar curvature, whereas the LeBrun metrics are (in suitable
domains) conformal to Pedersen metrics~\cite{Ped:emm}, which are AH
self-dual Einstein metrics of negative scalar curvature
(see~\cite{Hit:tem,CaPe:emt}).

Our main results show that a similar picture continues to hold for more
general cyclic subgroups of $U(2)$. We begin with the SFK metrics.

\begin{result} \label{thma} Let $\Gamma \subset U(2)$ be a finite cyclic
subgroup such that $0\in \C^2$ is the only fixed point, and let $X$ be a
toric resolution of $\C^2/\Gamma$.  Then $X$ admits a finite-dimensional
family of ALE scalar-flat K\"ahler metrics.
\end{result}

Several remarks are in order. First, this result extends to the case
$\Gamma=\{1\}$, when $X$ is an iterated blow-up of $\C^2$ and the metric is
asymptotically euclidean---in this form, the result is due to
Joyce~\cite{Joy:esd}. The above theorem follows easily from Joyce's work
and is also implicit in~\cite{Joy:hqqq,Joy:qcs}.

Second, these metrics are given by {\em explicit formulae} and have a
$2$-torus $T^2\cong S^1\times S^1$ acting by isometries. Note that there is
a minimal toric resolution $X_0$ of $\C^2/\Gamma$ (containing no $-1$
curves) such that any toric resolution of $\C^2/\Gamma$ is an iterated
blow-up of $X_0$ at fixed points of the torus action.  By the gluing
theorems of Kovalev and the second author~\cite{KoSi:gt}, there exist
complete scalar-flat K\"ahler metrics on {\em any} iterated blow-up of
$X_0$, but explicitness is then lost.

In our next result we find closely related non-K\"ahler SDE metrics, with
negative scalar curvature, which may be viewed as hyperbolic analogues of
the $A_n$ hyperk\"ahler gravitational instantons.

\begin{result} \label{thmb} Let $\Gamma$ and $X$ be as in
Theorem~\ref{thma} and suppose that $c_1(X) <0$.  Then there is a connected
open neighbourhood $X_+$ of the exceptional fibre $E$ of $X\to \C^2/\Gamma$
given by a smooth positive defining function $F$, with the property that
for {\em one} of the metrics $g$ in Theorem~\ref{thma},
\begin{equation} \label{eq:conffact}
g_+ := F^{-2}g
\end{equation}
is a complete self-dual Einstein metric with negative scalar curvature.
\end{result}

Here the zero-set
\begin{equation} \label{1.13.3}
Y = \{F=0\}
\end{equation}
is diffeomorphic to the link of the singularity at the origin in
$\C^2/\Gamma$, which is a lens space, and $X_+$ is diffeomorphic to $X$.
We remark that the constraint $c_1(X)<0$ does not give a restriction on the
size of $b_2(X)$, so this theorem supplies explicit complete Einstein
metrics on $4$-manifolds with arbitrarily large second Betti number.

\ifaddpics
\begin{figure}[ht]
\begin{center}
\includegraphics[width=.3\textwidth]{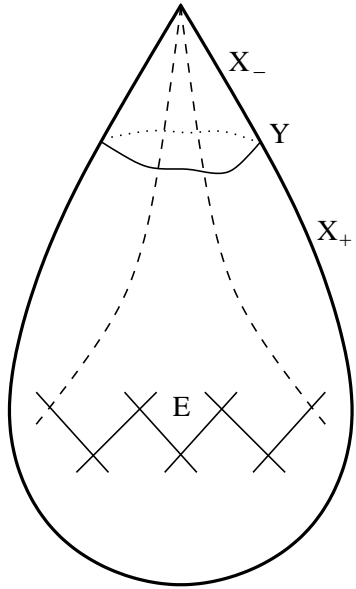}
\caption{The one point compactification of $X$.}
\label{fig0}
\end{center}
\end{figure}
\fi

The given formulation of Theorem~\ref{thmb} is designed to fit naturally
with the perspective of {\em asymptotically hyperbolic} metrics.  This is a
class of complete, conformally compact metrics generalizing the classical
relation
\begin{equation*}
g_{\mathrm{hyp}} = \frac{4}{(1-r^2)^2}g_{\mathrm{euc}}
\end{equation*}
between the hyperbolic metric $g_{\mathrm{hyp}}$ on the open ball $\{r<1\}$
and the euclidean metric $g_{\mathrm{euc}}$ on its closure. The key point
here is that after multiplication by the square of the defining function
$(1-r^2)/2$, $g_\mathrm{hyp}$ extends to a riemannian metric on the
boundary $\{r=1\}$.  The relation \eqref{eq:conffact} is an example of the
same kind, for $g$ extends smoothly (as a riemannian metric) to $X$.

More generally, if $M$ is a compact manifold with boundary $N$, we say that
a riemannian metric $g$ in the interior $M^o$ of $M$ is {\em asymptotically
\textup(locally\textup) hyperbolic} (AH) with {\em conformal infinity}
$(N,c)$, if for any boundary-defining function $u$, $u^2 g$ extends
smoothly to $N$ and $u^2 g$ is in the conformal class $c$.  In this
situation, we also say that $(M,g)$ is a {\em filling} of $(N,c)$ or that
$(N,c)$ {\em bounds} $(M,g)$. Notice that the freedom to multiply defining
functions by any function smooth and positive near $N$ is absorbed
precisely by the specification of a conformal class rather than a metric on
$N$.  These ideas suggest natural boundary-value problems, such as: given
$(N,c)$, does there exist a filling $(M,g)$ with $g$ and Einstein metric?
If so, is $g$ unique?  Following work of Fefferman--Graham, Graham--Lee,
Biquard and Anderson~\cite{FeGr:ci,GrLe:emc,Biq:mes,And:emc}, one is
beginning to have a good understanding of this problem: at least if $c$ has
positive Yamabe constant, it seems that `generically' the Einstein filling
exists and is unique up to diffeomorphism, so the problem is well-posed.

If $N$ is of dimension $3$, there is another boundary-value problem, namely
that of filling $c$ by a self-dual or anti-self-dual Einstein metric.  In
so far as the Einstein problem is well-posed, this problem must be
over-determined. An important step in the study of this problem is
Biquard's recent proof~\cite{Biq:mab} of the positive-frequency conjecture
of LeBrun~\cite{LeBr:hcc,LeBr:cqk}. This asserts (roughly speaking) the
existence of a decomposition $c = c_-+c_0+c_-$ (if $c$ is close to the
round conformal structure $c_0$ on $S^3$) with the property that $c_++c_0$
bounds a SDE metric on the ball and $c_-+c_0$ bounds an ASDE metric on the
ball---equivalently $c_++c_0$ and $c_-+c_0$ bound SDE metrics on the ball
inducing opposite orientations on $S^3$. This is a nonlinear version of the
decomposition of a function on the circle into positive and negative
Fourier modes.

The metric $g_+$ of Theorem~\ref{thmb} shows a behaviour similar to that of
the standard round metric on $S^3$: if $(Y,c)$ is the conformal infinity of
$g_+$, then $c$ bounds the SDE metric $g_+$ but also an SDE orbifold metric
$g_-$ (which is a cyclic quotient of a smooth metric on the ball) living on
the domain $X_- =\{F<0\}$ (the `other side' of $Y$).  This picture was well
known in the special case of the Pedersen--LeBrun metrics when $X$ is the
total space of $\cO(-p)$ (see e.g.~\cite{Hit:tem}).  It seems that there
are many new phenomena worthy of exploration here, to which we hope to
return.

Our next result belongs to this circle of ideas. To state it, we recall
that the link of the singularity in $X$ is a lens space.

\begin{result} \label{thmc}  Let $X$ be as above, with $c_1(X)<0$,
and let $N$ be the corresponding lens space. Then there is an
infinite-dimensional family of conformal structures on $N$ that bound
complete SDE metrics on neighbourhoods $X_+$ of the exceptional fibre $E$
of $X$.
\end{result}
These metrics are not quite as explicit as the one in Theorem~\ref{thmb},
since the required functions are in general given by integral formulae over
the boundary of hyperbolic space. They are parameterized (roughly speaking)
by distributions in one variable with compact support in $(-\infty,0)$.
Note also that the generic conformal structure in Theorem~\ref{thmc} will
{\em not} bound an SDE orbifold metric on $X_-$. Indeed, we show that there
are examples where the conformal structure on $X_+$ does not extend, as a
self-dual conformal structure, at all into $X_-$.

\vspace{5pt} The AH condition is not the only boundary condition we might
consider: Biquard~\cite{Biq:mes} also considers metrics which are
asymptotic to the Bergman metric of complex hyperbolic space. This may be
viewed as a limiting case in which the conformal structure on $N$
degenerates into a (pseudoconvex) CR structure.  Recall that the latter is
is given by a contact distribution $\cH\subset TY$ and a conformal
structure $c$ on $\cH$ which is compatible with the Levi form. This means
that there is a (uniquely determined) almost complex structure $J$ on
$\cH$, such that for any contact $1$-form $\theta$ (for $\cH$),
$d\theta\restr{\cH}$ is the K\"ahler form, with respect to $J$, of a metric
$h$ in $c$.

Now if $M$ is a compact manifold with boundary $N$, we say that a
riemannian metric $g$ in the interior $M^o$ of $M$ is {\em asymptotically
\textup(locally\textup) complex hyperbolic} (ACH) with {\em CR infinity}
$(N,\cH,c)$ if for any boundary-defining function $u$, there is a $1$-form
$\theta$ on $M$ such that $u^2 g-u^{-2}\theta^2$ extends to a smooth and
degenerate metric $h$ on $M$, and such that the pullback of $(\theta,h)$ to
$N$ is a contact metric structure compatible with $(\cH,c)$ on $N$. Note
again that the CR structure is independent of the choice of boundary
defining function (and the choice of the $1$-form $\theta$ on $M$).

The ACH boundary-value problem for Einstein metrics is studied by
Biquard~\cite{Biq:mes}.  As with the AH boundary condition, when $\dim
N=3$, we can strengthen the Einstein equation to the self-dual Einstein
equation.  To describe the explicit examples we obtain, we
recall~\cite{Bel:ncr} that a CR manifold $N$ is said to be {\it normal} if
it admits a Reeb vector field which generates CR automorphisms of $N$, and
in addition {\it quasiregular} if the orbit space is an orbifold Riemann
surface.

\begin{result}\label{thmch} Let $X$ be as above, with $c_1(X)<0$ and
let $N$ be the corresponding lens space. Then there is an ACH SDE metric
on $X$ whose CR infinity is a quasiregular normal CR structure on $N$.
\end{result}
In contrast to the AH boundary condition, we obtain only one SDE metric in
each case, \ie, we do not find deformations of the CR structure on $N$
which still bound ACH SDE metrics on manifolds diffeomorphic to $X$.

\vspace{5pt} Our results provide an intriguing analogue of the
K\"ahler--Einstein trichotomy: a compact K\"ahler manifold $M$ can only
admit a K\"ahler--Einstein metric if the first Chern class is
positive-definite, zero, or negative-definite, the sign of the scalar
curvature being respectively positive, zero or negative.  (Moreover, thanks
to the well-known work of Aubin, Calabi and Yau, the necessary conditions
are also sufficient, if $c_1(M)\leq 0$.)

We have seen here that our non-compact complex surface $X$ admits a toric
SDE metric of negative scalar curvature if $c_1(X)<0$ and admits a SDE
metric of zero scalar curvature if $c_1(X)=0$. Of course the only compact
SDE metrics of positive scalar curvature are $\CP2$ and $S^4$. Since the
negative-scalar-curvature SDE metrics are definitely not K\"ahler, it is
natural to ask whether there is an underlying reason for the condition
$c_1(X)<0$: for example, is this a special feature of the toric symmetry of
our metrics, or is there an analogous statement for general SDE metrics on
these surfaces?

Although we cannot answer this question in general, we can provide a
converse to Theorem~\ref{thmc} in the toric case, thanks to the following
result.

\begin{result} \label{thmd} Let $g$ be a toric SDE metric of
negative scalar curvature, defined in a connected neighbourhood $U$ of an
embedded sphere $S$ consisting of special orbits of the torus action.  Then
$|S\cdot S|\geq 3$.
\end{result}

In particular, if $X$ is the complex manifold underlying one of the $A_n$
hyperk\"ahler instantons, then $X$ does not admit a toric SDE metric of
negative scalar curvature.

\vspace{5pt} The methods of this paper are essentially explicit and depend
crucially on the classification of toric SDE metrics of nonzero scalar
curvature by Pedersen and the first author~\cite{CaPe:emt}. It is indeed
remarkable that the torus-symmetry reduces the SDE equations to a standard
linear partial differential equation in the hyperbolic plane, which can be
studied relatively easily. This may be contrasted with the case of $SU(2)$
or $SO(3)$ symmetry which leads to a nonlinear ordinary differential
equation~\cite{Hit:tem,Tod:p6}. Solutions of the latter generally have to
be written down in terms of $\vartheta$-functions.

\vspace{5pt} In the first two sections we review some essential background:
the geometry of resolutions of cyclic singularities from the complex and
smooth points of view; and the local constructions of $T^2$-invariant
self-dual and SDE metrics. Then we move on to the proofs of the main
theorems giving some examples along the way. An appendix is devoted to a
brief exposition of some aspects of the geometry of hyperbolic space that
were suppressed in the body of the paper.

\subsection*{A note on orientations} If $X$ is a complex surface, then
with the standard complex orientation, a scalar-flat K\"ahler metric is
anti-self-dual. The Fubini--Study metric on $\CP{2}$, again with the
standard complex orientation, is self-dual (and Einstein). In this paper
are primary concern is (anti-)self-dual Einstein metrics, which generalize
the Fubini--Study metric in a natural sense. We have therefore decided to
state our results for self-dual metrics, which implies a reversal of the
complex orientation of the complex surface $X$ in the above. We hope the
reader will not find this change of orientations too tiresome.

\acknowledge Thanks to Henrik Pedersen, Kris Galicki, Olivier Biquard,
Michael Anderson and Rafe Mazzeo for helpful discussions. The first author
is grateful to the Leverhulme Trust and the William Gordan Seggie Brown
Trust for financial support. Both authors are members of EDGE, Research
Training Network HPRN-CT-2000-00101, supported by the European Human
Potential Programme

The diagrams were produced using Xfig, MAPLE and Paul Taylor's
commutative diagrams package.
%
% Section 2: Background on Hirzebruch-Jung strings and the corresponding 
% singularities and the combinatorial description of toric 4-manifolds
%
\section{Toric topology}

\subsection{Hirzebruch-Jung strings and cyclic surface singularities}

Let $p$ and $q$ be coprime integers with $p > q>0$. Let $\Gamma$ be
the cyclic subgroup of $U(2)$ generated by the matrix
\begin{equation} \label{21.12.1}
\begin{bmatrix} \exp(2\pi i/p) & 0 \cr
           0 & \exp(2\pi iq/p)\cr\end{bmatrix}.
\end{equation}
The quotient $\C^2/\Gamma$ is a complex orbifold, with an isolated singular
point corresponding to the origin. The resolution of this singularity by a
Hirzebruch--Jung string is well known in algebraic geometry, cf.~\cite[\S
III.5]{BPV:ccs}, especially \cite[Theorem~5.1 and
Proposition~5.3]{BPV:ccs}. In outline, the story is as follows.  There is a
{\em minimal}\footnote{\ie, any other resolution is some blow-up of this
one, or equivalently, there are no $(-1)$-curves} resolution $\pi:X \to
\C^2/\Gamma$ with the properties
\begin{enumerate}
\item $X$ is a smooth complex surface;
\item there is an {\em exceptional divisor} $E\subset X$ such that 
$\pi(E) = \{0\}$, but the restriction of $\pi$ is a biholomorphic map
$X\setdif E \to (\C^2/\Gamma) \setdif \{0\}$;
\item $E =S_1\cup S_2\cup\cdots\cup S_k$ where the $S_j$ are 
holomorphically embedded smooth $2$-spheres and the intersection
matrix of the $S_j$ has the form
\begin{equation} \label{21.12.2}
(S_i\cdot S_j) = \begin{bmatrix} -e_1 & 1 & 0 &  \cdots & 0 \cr 1 &
-e_2 & 1 & \cdots & 0 \cr 0 & 1 &-e_3 &  \cdots & 0 \cr \vdots
& \vdots & \vdots &  & \vdots \cr 0 & 0 & 0 & \cdots &
-e_k \cr \end{bmatrix}
\end{equation}
where all $e_j \geq 2$.
\end{enumerate}
The integers $e_j$ are determined by $p$ and $q$ through the
continued-fraction expansion
\begin{equation} \label{21.12.3}
\frac{q}{p} = \frac{1}{e_1-}\frac{1}{e_2-}\cdots\frac{1}{e_k}
\end{equation}
or equivalently by the euclidean algorithm in the form:
\begin{equation} \label{21.12.4}
r_{-1}:=p,\qquad r_0:=q,\qquad r_{j-1} = r_j e_{j+1} - r_{j+1},\quad
\mbox{where}\quad  0 \leq r_{j+1} < r_j.
\end{equation}
Note that in this version of the division algorithm, the quotients $e_j$
are overestimated and since $r_{j+1} <r_{j}$, we have $e_j\geq 2$ for
$j=1,\ldots k$.  We remark also that the effect of reversing the order of
the $e_j$ in this continued fraction is to replace $q/p$ by $\tilde{q}/p$,
where $q\tilde{q} \equiv 1$ mod $p$. This replacement in the action
\eqref{21.12.1} does not change the orbifold singularity, just as one would
expect.

By the adjunction formula, $c_1(X)\leq 0$; we have $c_1(X)<0$ if all
$e_j>2$ and $c_1(X)=0$ if all $e_j=2$. In the latter case, $q = p-1$,
$\Gamma\subset SU(2)$ and $X$ is the canonical resolution of the $A_{p-1}$
singularity.

\subsection{Toric differential topology}\label{sec:toric}
We turn now to a discussion of `toric' $4$-manifolds and orbifolds. These
are a class of smooth $4$-manifolds (orbifolds) with a given smooth action
of the $2$-torus $T^2= S^1\times S^1 = \R^2/2\pi\Z^2$ with the property
that the generic orbit is a copy of $T^2$. The source for this material is
the paper of Orlik and Raymond \cite{OrRa:at4}; the summary in
\cite{Joy:esd} is also useful.

We begin with the standard action of $T^2$ on $\C^2=\R^4$
\begin{equation} \label{1.24.2}
(t_1,t_2)\cdot(z_1,z_2) = (e^{it_1}z_1,e^{it_2}z_2)
\end{equation}
where $(t_1,t_2)$ are standard coordinates on $T^2$. The quotient space
$\R^4/T^2$ is easily described by introducing polar coordinates $z_j =
r_je^{i\theta_j}$. Then it is clear that $\R^4/T^2$ is the quarter-space
$Q=\{(r_1,r_2):r_1\geq 0,r_2\geq 0\}$ and that the orbit corresponding to
any interior point of $Q$ is a copy of $T^2$, while the orbits
corresponding to boundary points have non-trivial stabilizers inside $T^2$.
Since the origin is the only fixed point of this $T^2$-action, the
stabilizer of the corner $(0,0)\in Q$ is the whole of $T^2$, while the
stabilizer of orbits corresponding to the positive $r_1$ and $r_2$ axes are
certain circle subgroups of $T^2$. Such circle subgroups are essential in
what follows and we adopt the description of them used by Orlik and
Raymond.

\subsubsection{Definition} For any pair of coprime integers $(m,n)$,
$G(m,n) =\{mt_1+nt_2=0\}\subset T^2$.
(It is clear that $G(m,n) = G(\tilde m,\tilde n)$ if and only if 
$(\tilde m,\tilde n) = \pm(m,n)$.)

\vspace{5pt} With these conventions the stabilizer of the positive
$r_2$-axis in $\R^4$ is $G(1,0)$ and that of the positive $r_1$-axis is
$G(0,1)$. A complete description of the quotient space consists of the
quarter space, with its two edges labelled by their stabilizers:
\begin{diagram}[silent,width=2.2em,abut]
\null&\isot{(1,0)}{\tilde S}&\dotify{\mathop\bullet\limits_{}^{}}
&\isot{(0,1)}{S}&\null
\end{diagram}
Here we have opened out the quarter-space in order to simplify the diagram
and have labelled the edges (which correspond to the special orbits with
circle-stabilizers) $S$ and $\tilde S$.  It is to be understood with such
diagrams that the generic orbits correspond to the upper half-space.

This $T^2$-action on $\R^4$ extends smoothly to the one-point
compactification $S^4$, the point at infinity being another fixed point.
The quotient $S^4/T^2$ is obtained as the one-point compactification of
$Q$---it is a di-gon with two vertices and two edges labelled by $(1,0)$ and
$(0,1)$.

\vspace{5pt} Now consider the orbifold $\R^4/\Gamma$. Since the action of
the finite group $\Gamma$ in \eqref{21.12.1} commutes with the action of
$T^2$ on $\R^4$ it follows that $\R^4/\Gamma$ is a toric orbifold.  The
quotient map $\R^4 \to \R^4/\Gamma$ induces a map $Q\times T^2 \to Q\times
T^2$, which must be the identity on $Q$. The required map $T^2\to T^2$ must
be surjective with kernel $\Z(1/p,q/p)$ so can be taken as
\begin{equation*}
(t_1,t_2) \mapsto (pt_1, qt_1-t_2).
\end{equation*}
The labels $(1,0)$ and $(0,1)$ map to $(p,q)$ and $(0,-1)$ by this map, so
the combinatorial picture of $\R^4/\Gamma$ is as follows:
\begin{diagram}[silent,width=2.2em,abut]
\null&\isot{(p,q)}{\tilde S}&\dotify{\mathop\bullet\limits_{}^{}}
&\isot{(0,-1)}{S}&\null
\end{diagram}
The compactification $S^4/\Gamma$ is an orbifold with two isolated singular
points and the corresponding quotient is a di-gon with edges labelled by
$(p,q)$ and $(0,-1)$.

\vspace{5pt} According to Orlik and Raymond, the general picture of a
smooth $4$-orbifold $M$ with isolated singular points and smooth
$T^2$-action will be modelled by the examples we have just described. More
precisely, if we assume that $M$ is simply connected, then the quotient
$M/T^2$ is topologically a closed polygon (simply connected 2-manifold with
corners). A typical edge $S_j$ is labelled by the coprime pair $(m_j,n_j)$
such that $G(m_j,n_j)$ is the corresponding isotropy group. It is
convenient also to label the vertex $S_j\cap S_{j+1}$ by the number
\begin{equation} \label{2.24.2}
\eps_j = m_jn_{j+1} -m_{j+1}n_j
\end{equation}
so that the full picture of the boundary of $M/T^2$ is as follows:
\begin{diagram}[silent,width=2.2em,abut]
\cdots
\null&
\isot{\;\;(m_{j+1},n_{j+1})\;\;\;\;}{S_{j+1}}&
\dotify{\mathop\bullet\limits_{}^{\ve_j}}&
\isot{\;\;\;(m_{j},n_{j})\;\;\;}{S_{j}}&
\dotify{\mathop\bullet\limits_{}^{\ve_{j- 1}}}&
\isot{\;\;\;\;(m_{j-1},n_{j-1})\;\;}{S_{j-1}}&
\null\cdots
\end{diagram}
The union of orbits corresponding to the edge $S_j$ is an embedded
$2$-sphere in $M$ and the vertex $S_j\cap S_{j+1}$ is a smooth point of $M$
if $\ve_j= \pm 1$ and is more generally an orbifold point with isotropy of
order $|\ve_j|$.  Indeed the precise quotient singularity can be discovered
by transforming $(m_j,n_j)$ to $(0,-1)$ and $(m_{j+1}, n_{j+1})$ to
$(p,q)$, $p>q>0$, by an element of $SL(2,\Z)$, and comparing with the model
described above. One of the crucial points of the work of Orlik-Raymond is
the converse to the above description, namely that {\em any} diagram of the
above type gives rise to a unique $4$-orbifold with a given action of
$T^2$.

The topological properties of $M$ are encoded in diagrams of the above
kind. For example, let
\begin{equation}\label{2.13.3}
e_j :=\ve_{j-1}\ve_{j}(m_{j-1}n_{j+1} - m_{j+1}n_{j-1}).
\end{equation}
Then if $|\ve_{j-1}| = |\ve_{j}| = 1$ (and with a suitable orientation
convention),
\begin{equation} \label{1.19.3}
 S_j\cdot S_j = e_j.
\end{equation}
Notice that here we have blurred the distinction between the edge $S_j$ and
the corresponding sphere in $M$ which lies over this edge. We shall
continue with this abuse wherever convenient.

\subsection{Toric resolutions of singularities}
In order to resolve the singularity at the origin of $\R^4/\Gamma$
\begin{diagram}[silent,width=2.2em,abut]
\null&\isot{(p,q)}{\tilde S}&\dotify{\mathop\bullet\limits_{}^{p}}
&\isot{(0,-1)}{S}&\null
\end{diagram}
we must replace the vertex by a chain of vertices and edges of the
following kind
\begin{diagram}[silent,width=2.2em,abut]
\null
&\isot{(p,\tilde q)}{\tilde S=S_{k+1}}&\dotify{\mathop\bullet\limits_{}^{1}}
&\isot{(m_k,n_k)}{S_{k}}&\dotify{\mathop\bullet\limits_{}^{1}}
&\rLine&{\,\cdots\,} &\rLine &\dotify{\mathop\bullet\limits_{}^{1}}
&\isot{(m_2,n_2)}{S_2}&\dotify{\mathop\bullet\limits_{}^{1}}
&\isot{(1,0)}{S_{1}}&\dotify{\mathop\bullet\limits_{}^{1}}
&\isot{(0,-1)}{S=S_0}&\null
\end{diagram}
where $\tilde q \equiv q \pmod{p}$ and we have used up the freedom to change
bases by choosing
\begin{equation*}
(m_0,n_0) = (0,-1),\quad (m_1,n_1) = (1,0).
\end{equation*}
(This normalization, which will be convenient later, forces us to allow $q$
to be replaced by $\tilde q$ as above.  We have also chosen all
$\eps_j=+1$, which can always be done by playing with the sign ambiguity of
the $(m_j,n_j)$.)

For the rest of this paper we shall assume also that $m_j>0$ for
$j=1,\ldots k+1$. This has the effect of making the toric resolution have
semi-definite intersection-form, as follows from \eqref{1.19.3}. Thus we
make the following definition.
\subsubsection{Definition} A sequence of coprime integers $(m_j,n_j)$
is called {\em admissible} if
\begin{enumerate}
\item $(m_0,n_0) = (0,-1)$, $(m_1,n_1) = (1,0)$,
$(m_{k+1},n_{k+1})=(p,\tilde q)$, $\tilde q \equiv q \pmod{p}$; 
\item $m_j>0$ for $j=1,2,\ldots k+1$;
\item $m_jn_{j+1}-m_{j+1}n_j =1$ for $j=0,1,\ldots  k$.
\end{enumerate}
Then there is a one to one correspondence between admissible sequences and
toric resolutions of $\R^4/\Gamma$ with semi-definite intersection form.

\vspace{5pt} We note also that the compactified version of this resolution
is given by the diagram
\begin{diagram}[silent,width=2.2em,abut]
\null
&\isot{(0,1)}{S_{0}}&\dotify{\mathop\bullet\limits_{}^{p}}
&\isot{(p,\tilde q)}{S_{k+1}}&\dotify{\mathop\bullet\limits_{}^{1}}
&\isot{(m_k,n_k)}{S_{k}}&\dotify{\mathop\bullet\limits_{}^{1}}
&\rLine&{\,\cdots\,}
&\rLine
&\dotify{\mathop\bullet\limits_{}^{1}}
&\isot{(m_2,n_2)}{S_2}&\dotify{\mathop\bullet\limits_{}^{1}}
&\isot{(1,0)}{S_{1}}&\dotify{\mathop\bullet\limits_{}^{1}}
&\isot{(0,-1)}{S_0}&\null
\end{diagram}
where the point labelled $p$ is the orbifold point at infinity.

We observe that for an admissible sequence,
\begin{equation*}
\frac{n_{j+1}}
{m_{j+1}} = \frac{1}{e_1-}\frac{1}{e_2-}
\cdots\frac{1}{e_j}\quad\mbox{for}\quad j=1,2,\ldots k,
\end{equation*}
where the $e_j$ are defined by \eqref{2.13.3}. Conversely, the
continued-fraction expansion \eqref{21.12.3} defines the {\em minimal
admissible sequence} and this is the toric description of the
Hirzebruch-Jung string.  It is not hard to check that any admissible
sequence arises from the minimal one by blow-up, (\ie, by connected sum
with $\CP{2}$) which corresponds to the replacement of
\begin{diagram}[silent,width=2.2em,abut]
\cdots
\null&
\isot{\;\;(m_{j+1},n_{j+1})\;\;\;\;}{S_{j+1}}&
\dotify{\mathop\bullet\limits_{}^{\pm 1}}&
\isot{\;\;\;(m_{j},n_{j})\;\;\;}{S_{j}}&
\null\cdots
\end{diagram}
by
\begin{diagram}[silent,width=2.2em,abut]
\cdots
\null&
\isot{\;\;(m_{j+1},n_{j+1})\;\;\;\;}{S_{j+1}}
&\dotify{\mathop\bullet\limits_{}^{\pm
1}}&
\isot{\;\;\;(m_{j+1}+m_{j},n_{j+1}+n_j)\;\;\;}{E}&
\dotify{\mathop\bullet\limits_{}^{\pm 1}}&
\isot{\;\;\;\;(m_{j},n_{j})\;\;}{S_{j}}&
\null\cdots
\end{diagram}
where we have labelled the exceptional divisor $E$ and have abused notation
by labelling the proper transforms of $S_j$ and $S_{j+1}$ by the same
symbols.  In terms of continued fractions, this process of blowing up
corresponds to the identity
\begin{equation*}
\cdots \frac{1}{a-}\;\frac{1}{b-}\cdots =
\cdots \frac{1}{a+1-}\;\frac{1}{1-}\;\frac{1}{b+1-}\cdots.
\end{equation*}
If a blow-up is performed at either end of the sequence, then $q$ will be
changed by a multiple of $p$.

%
% Section 3: Background on spin geometry of the hyperbolic plane, and the
% explicit classification of toric SDE metrics.
%
\section{Local constructions of half-flat conformal structures}

In this section we review the local construction of half-flat conformal
structures admitting two commuting Killing vectors. For more details,
please see \cite{CaPe:emt,Joy:esd}.

\subsection{Joyce's construction}

\subsubsection{Data}
\begin{itemize}
\item $(B,h)$ is a spin $2$-manifold with metric $h$ of constant 
curvature $-1$;

\item $\cW \to B$ is the spin-bundle, viewed as a real vector bundle of
rank $2$ equipped with the induced metric, also denoted $h$;

\item $\V$ is a given $2$-dimensional vector space, with symplectic
structure $\ip{\cdot,\cdot}$ and $M = B\times \V$ is the corresponding
product bundle.
\end{itemize}
Note that $\V$ acts in the obvious way by translations on $M$.

\vspace{5pt} Given data as above, we define a $\V$-invariant riemannian
metric on $M$ associated to any bundle isomorphism $\Phi\colon \cW \to
B\times \V$.  Indeed, given $\Phi$, we define a family of metrics on $\V$,
parameterized by $B$,
\begin{equation*}
(v,\tilde v)_\Phi = h(\Phi^{-1}(v),\Phi^{-1}(\tilde v))
\end{equation*}
and then a metric on the total space
\begin{equation} \label{27.1.2}
g_\Phi = \Omega^2(h + (\cdot,\cdot)_\Phi),
\end{equation}
where the conformal factor $\Omega>0$ is in $\Cinf(B)$. It is clear that
such a metric is invariant under the action of $\V$ on $M$. Any such metric
also descends to the quotient $B\times (\V/\Lambda)$ where $\Lambda\subset
\V$ is any lattice.

The remarkable observation of Joyce is that $g_\Phi$ is conformally
half-flat if $\Phi$ satisfies a linear differential equation that we shall
call here the Joyce equation. The Joyce equation essentially makes $\Phi$
an eigenfunction of the Dirac operator. Here we shall content ourselves
with an explicit form of this equation, suitable for our later
purposes. (See the appendix for a sketch of the underlying geometry.)

\subsubsection{The Joyce equation} If $B\subset \cH^2$ we can 
introduce half-space coordinates $(\rho,\eta)$, with $\rho>0$, so that
\begin{equation*}
h =\frac{d\rho^2 + d\eta^2}{\rho^2}
\end{equation*}
and an orthonormal frame $\lam_1$, $\lam_2$ of $\cW^*$ such that
\begin{equation*}
\lam_1\otimes\lam_1 - \lam_2\otimes\lam_2
= d\rho/\rho,\quad \lam_1\otimes\lam_2 + \lam_2\otimes\lam_1 =  d\eta/\rho.
\end{equation*}
(From the complex point of view, $\lam_1+i\lam_2 =
\sqrt{(d\rho+id\eta)/\rho}$.)

Since $\Phi \in \Cinf(B,\cW^*\otimes \V)$ we can write
\begin{equation*}
\Phi = \lam_1\otimes v_1 + \lam_2 \otimes v_2
\end{equation*}
where $v_1$ and $v_2$ are in $\Cinf(B,\V)$. Then the Joyce equation 
is the system
\begin{equation} \label{27.1.6}
\rho\del_\rho v_1 + \rho\del_\eta v_2 = v_1,\quad
\rho\del_\eta v_1 - \rho\del_\rho v_2 =0.
\end{equation}
Clearly $\Phi$ defines a bundle isomorphism if $\ip{v_1,v_2}\neq0$ and
then
\begin{equation*}
\Phi^{-1}=\frac{\ip{v_1,\cdot}\tens\ell_1-\ip{v_2,\cdot}\tens\ell_2}
{\ip{v_1,v_2}},
\end{equation*}
where $\ell_1,\ell_2$ is the dual orthonormal frame of $\cW$.

Let us summarize these observations.

\begin{thm}\noindent~\tcite{Joy:esd} With the above notation, if 
\begin{equation} \label{27.1.9}
\ip{v_1,v_2} \neq 0\mbox{ in }B
\end{equation}
and $v_1$ and $v_2$ satisfy \eqref{27.1.6}, then for any conformal 
factor $\Omega$, the metric
\begin{equation} \label{27.1.7}
g_\Phi = \Omega^2\Bigl(
\frac{d\rho^2 + d\eta^2}{\rho^2}
+ \frac{\ip{v_1,\cdot}^2 + \ip{v_2,\cdot}^2}{\ip{v_1,v_2}^2}
\Bigr)
\end{equation}
is a $\V$-invariant conformally half-flat metric on $M$.
\end{thm}

\subsubsection{Scalar-flat K\"ahler representatives}

For each choice of a point $y$ on the circle at infinity of $\cH^2$, there
is a scalar-flat K\"ahler metric conformal to $g_\Phi$:
\begin{equation} \label{17.1.2}
g_{SFK}=\frac{\rho|\ip{v_1,v_2}|}{\rho^2+(\eta-y)^2}
\Bigl(\frac{d\rho^2+d\eta^2}{\rho^2}
+ \frac{\ip{v_1,\cdot}^2 + \ip{v_2,\cdot}^2}{\ip{v_1,v_2}^2}
\Bigr).
\end{equation}
The fact that $g_{SFK}$ is a scalar-flat K\"ahler metric follows
from~\cite[Proposition 2.4.4]{Joy:esd}, or alternatively, as explained
in~\cite{CaPe:emt}, from the work of LeBrun~\cite{LeBr:cp2}.

\subsubsection{Basic solutions of the Joyce equation}

There is a basic solution of the Joyce equation \eqref{27.1.6}---or rather
of the corresponding equation where the coefficients $v_1$ and $v_2$ are
ordinary functions rather than $\V$-valued functions---associated to any
given point $(0,y)$ of the boundary of $\cH^2$:
\begin{equation} \label{27.1.8}
\phi(\rho,\eta;y)
= \frac{\rho\lam_1 + (\eta -y)\lam_2}{\sqrt{\rho^2 + (\eta-y)^2}}.
\end{equation}
In \cite{Joy:esd} solutions of \eqref{27.1.6} given by finite linear
combinations of the form $\phi(\rho,\eta;y)\otimes v$ were used to
construct conformally half-flat metrics on connected sums of the complex
projective plane. This idea will be used in this paper too, though we shall
also use suitable infinite linear combinations of these solutions to
construct infinite-dimensional families of metrics.

\subsubsection{Compactification}\label{sec:cpt}

We now explain how this construction of conformally half-flat metrics is
combined with the previous description of toric $4$-manifolds. The main
point is to give sufficient conditions for the smooth extension of the
metric $g_\Phi$ over the special orbits at the boundary.  The following
sufficient conditions were given by Joyce. Consider a combinatorial diagram
with labelled edges and vertices as in \S\ref{sec:toric} and consider the
boundary of $M/T^2$ to be identified with the boundary $\rho=0$ of $\cH^2$.
Note that $\V$ is the Lie algebra of $T^2$, which we identify with $\R^2$.
Then if $B$ is a neighbourhood in $\cH^2$ of a point in the interior of an
edge labelled $(m,n)$, the metric \eqref{27.1.7} extends smoothly to the
special orbits over $S\cap B$ provided that
\begin{equation} \label{27.1.10}
v_1 = O(\rho),  \quad v_2 = (m,n) + O(\rho^2)\mbox{ as }
\rho\to 0
\end{equation}
and 
\begin{equation} \label{27.1.11}
\rho^{-2}\Omega^2\mbox{ is smooth and positive in } B.
\end{equation}
Similarly, if $B$ is now a neighbourhood of a corner at $\eta=\eta_0$,
say, then \eqref{27.1.7} extends smoothly to the fixed-point if in
addition to the boundary conditions \eqref{27.1.10} and
\eqref{27.1.11}, the conformal factor is chosen so that
\begin{equation} \label{27.1.13}
\frac{\sqrt{\rho^2 + (\eta -\eta_0)^2}}{\rho^2}\Omega^2
\end{equation}
is smooth near the corner at $\eta=\eta_0$.

\subsubsection{Boundary behaviour of the basic solution}
The basic solution $\phi(\rho,\eta;y)$ has very simple boundary behaviour,
so that the preceding sufficient conditions are easily checked for linear
combinations of them:
\begin{equation} \label{28.1.1}
\mbox{if }\eta\neq y,\quad \phi(\rho,\eta;y) = O(\rho)\lam_1 + 
\bigl(\sign(\eta-y) + O(\rho^2)\bigr)\lam_2\mbox{ for small }\rho>0.
\end{equation}

\subsection{Self-dual Einstein metrics}
\label{sec:loc}
Return now to the local considerations. Suppose that $F\in \Cinf(B)$ 
is an eigenfunction with eigenvalue $3/4$ of the hyperbolic laplacian
\begin{equation} \label{27.1.15}
\Delta F = \frac{3}{4}F.
\end{equation}
Such a function is a {\em potential} for a Joyce matrix in the
following sense. Set
\begin{equation} \label{27.1.16}
f(\rho,\eta) = \sqrt{\rho}F(\rho,\eta),\quad
v_1 = (f_\rho,\eta f_\rho-\rho f_\eta),\quad
v_2 = ( f_\eta,\rho f_\rho +\eta f_\eta -f),
\end{equation}
regarded as $\V$-valued functions on $B$. Then, as is easily verified,
\begin{equation*} 
\Phi =\frac{1}{2}( \lam_1\otimes v_1 + \lam_2\otimes v_2)
\end{equation*}
is a solution of the Joyce equation.  The significance of these special
solutions of Joyce's equation is as follows.

\begin{thm}\label{thm:cp}~\tcite{CaPe:emt}
Let $F$ be as in \eqref{27.1.15} and $v_1$ and $v_2$ as in \eqref{27.1.16}.
Suppose further that
\begin{equation} \label{27.1.21}
F^2 \neq 4|d F|^2
\end{equation}
in $B$.  Let
\begin{equation*}
 D_+ = \{F>0\}, \quad Z= \{F=0\}, \quad D_- = \{F<0\}.
\end{equation*}
Then the metric
\begin{equation} \label{27.1.20}
g_F = \frac{\bigl|F^2-4|d F|^2\bigr|}{4F^2}
\left(
\frac{d\rho^2 + d\eta^2}{\rho^2}
+ \frac{\ip{v_1,\cdot}^2 + \ip{v_2,\cdot}^2}{\ip{v_1,v_2}^2}
\right)
\end{equation}
is a $\V$-invariant self-dual Einstein metric in $D_+\times\V$ and $D_-
\times \V$ \textup(if non-empty\textup). The scalar curvature has the same
sign as $F^2 - 4 |d F|^2$, and conversely any $\V$-invariant SDE metric of
nonzero scalar curvature arises \textup(locally, and up to
homothety\textup) in this way.
\end{thm}

Notice that $Z$ is non-empty then the scalar curvature must be negative and
that \eqref{27.1.21} ensures that $F$ is a defining function for $Z$. Since
$F^2g_F$ extends to a smooth metric in $B\times\V$, it is clear that $g_F$
is `locally conformally compact' with $Z\times \V$ as its conformal
infinity.

\subsubsection{Basic solutions}
For each $y$, the function
\begin{equation} \label{27.1.22}
F(\rho,\eta;y) = \frac{\sqrt{\rho^2 + (\eta-y)^2}}{\sqrt{\rho}}
\end{equation}
is an eigenfunction satisfying \eqref{27.1.15} and it is easy to check that
$F(\rho,\eta;y)$ is a potential for $\phi(\rho,\eta;y)\otimes(1,y)$.

\subsubsection{Remark}\label{sss1}
The role of $F$ as a potential for $\Phi$ relies on some geometry which we
outline in the appendix. The idea is that $\cW$ is being identified with
$\cH^2\times\V$; the derivative $dF$ of $F$ is then a section of
$\cW^*\tens\V^*$, which gives a section of $\cW^*\tens\V$ using the
symplectic form $\ip{\cdot,\cdot}$. Explicitly, as a section of $S^2\cW^*$,
$\Phi$ is given by
\begin{equation*}
\Phi = (\tfrac{1}{2}F + \rho F_\rho)\lam_1\tens\lam_1 +
(\rho F_\eta)(\lam_1\tens\lam_2+\lam_2\tens\lam_1) +
(\tfrac{1}{2}F - \rho F_\rho)\lam_2\tens\lam_2.
\end{equation*}
The identification of $\cW^*$ with $\cH^2\times\V$ is then obtained by
setting
\begin{equation*}
\lam_1 = (1/\sqrt{\rho},\eta/\sqrt{\rho}),\quad
\lam_2 = (0,-\sqrt{\rho}).
\end{equation*}
Note that $\det\Phi = \tfrac 14 F^2 - |dF|^2$ which is used
in~\eqref{27.1.21}.

\subsubsection{Example: hyperbolic space} The hyperbolic metric on the 
unit ball in $\R^4$ takes the form
\begin{equation*}
g = (1-r_1^2-r_2^2)^{-2}(dr_1^2 + dr_2^2 + r_1^2d\vartheta_1^2
+r_2^2d\vartheta_2^2)
\end{equation*}
in coordinates $(r_1,\vartheta_1)$, $(r_2,\vartheta_2)$ adapted to the 
standard action of $T^2$. This metric arises in the above construction 
from the function 
\begin{equation*}
F = \frac{1}{2}\bigl(F(\rho,\eta;-1)- F(\rho,\eta;1)\bigr).
\end{equation*}

To make this explicit is a matter of direct computation, using the 
formula
\begin{equation*}
(r_1+ir_2)^2 = \frac{\eta-1 +i\rho}{\eta+1 +i\rho}.
\end{equation*}
The reader may care to verify that
\begin{equation*}
\det\Phi = -\frac{\rho}{\sqrt{\rho^2+(\eta-1)^2}
\sqrt{\rho^2+(\eta+1)^2} }
\end{equation*}
and that the fibre metric (including the conformal factor) is
\begin{equation*}
\frac{\sqrt{\rho^2 +(\eta-1)^2}\sqrt{\rho^2 +(\eta+1)^2}}
{\bigl(\sqrt{\rho^2 +(\eta-1)^2}- \sqrt{\rho^2 +(\eta+1)^2}\,\bigr)^2}
\begin{bmatrix} (1 - u)/2 & 0\cr 0 & (1+u)/2
\end{bmatrix}
\end{equation*}
where
\begin{equation*}
u = \frac{\rho^2 +\eta^2-1}{
\sqrt{\rho^2 +(\eta-1)^2} \sqrt{\rho^2 +(\eta+1)^2}}.
\end{equation*}
Multiplying out and changing from $(\rho,\eta)$ to $(r_1,r_2)$ yields 
the hyperbolic metric.

%
%  The existence of canonical metrics: scalar-flat Kaehler on
%  resolutions of Hirzebruch-Jung singularities and the canonical SDE 
%  metric.
%
\section{The SFK metrics and canonical SDE metrics}

In the following three sections, we shall prove the theorems stated in the
introduction. In this section we shall give rather brief indications of the
proofs of Theorems~\ref{thma} and~\ref{thmb} (the theorems concerning the
scalar-flat K\"ahler metrics and canonical SDE metrics associated to
Hirzebruch--Jung strings).  Further details are given, in a more general
setting, in the next section. Indeed the results given here can be read as
a guide to the next section, as an extended set of examples.

\subsection{ALE scalar-flat K\"ahler metrics}

In order to prove Theorem~\ref{thma}, we take up the toric description of
Hirzebruch--Jung strings (and blow-ups of these strings) from
\S\ref{sec:toric}. Thus we suppose given an admissible sequence
$(m_j,n_j)$ for $j=0,\ldots k+1$, where
\begin{equation*}
(m_{k+1},n_{k+1}) =  (p,\tilde q), \quad \tilde q \equiv q \pmod{p}.
\end{equation*}
We supplement this by setting
\begin{equation*}
(m_{k+2},n_{k+2}) = -(m_0,n_0) = (0,1).
\end{equation*}
Since 
\begin{equation*}
m_jn_{j+1} - m_{j+1}n_j =
1\quad (j\neq k+1),\quad m_{k+1}n_{k+2}-m_{k+2}n_{k+1} = p
\end{equation*}
we see that the sequence of rationals $n_j/m_j$ enjoys a {\em
monotonicity property}
\begin{equation*} 
\frac{n_{j+1}}{m_{j+1}} > \frac{n_{j}}{m_{j}} \mbox{ all }j.
\end{equation*}

\subsubsection{Notation} Throughout this section we shall denote by $X$ the 
complex manifold corresponding to the data $(m,n)$---as we have seen, it is
a resolution of $\C^2/\Gamma$ obtained from the minimal one by blowing up
at fixed points of the $T^2$-action. We denote by $\oX$ the one-point
compactification (the added point will be denoted $\infty$).

\vspace{5pt}
Given an admissible sequence, pick real numbers
$y_0>y_1>\cdots>y_{k+1}$ and define
\begin{equation} \label{28.1.2}
\Phi(\rho,\eta)
= \frac{1}{2}\sum_{j=0}^{k+1} 
\phi(\rho,\eta;y_j)\otimes(m_j-m_{j+1},n_j-n_{j+1}).
\end{equation}
As a finite sum of the basic solutions $\phi(\rho,\eta;y)$, it is 
clear that $\Phi$ satisfies the Joyce equation and hence the 
metric $g_\Phi$ of \eqref{27.1.2} is a toric conformally half-flat metric 
wherever $\det\Phi\neq 0$. We resum \eqref{28.1.2} in the form
\begin{equation} \label{28.1.5}
\Phi(\rho,\eta) = \frac{1}{2}\sum_{j=0}^{k+1}\bigl(\phi(\rho,\eta;y_{j}) -
\phi(\rho,\eta;y_{j-1})\bigr)\otimes (m_j,n_j)
\end{equation}
where we decree
\begin{equation*}
\phi(\rho,\eta;y_{-1}) = - \phi(\rho,\eta;y_{k+1});
\end{equation*}
it then follows from \cite[Lemma 3.3.2]{Joy:esd} and the monotonicity of 
the rationals $n_j/m_j$ that $\det\Phi<0$ in $\cH^2$. (See also
Proposition~\ref{prop:neg} below.)
So $g_\Phi$ defines a $T^2$-invariant conformally 
half-flat metric over the whole of $T^2\times \cH^2$. Now it follows 
easily from the definitions and the known boundary behaviour of 
$\phi(\rho,\eta;y)$ that 
if $y_{j}<\eta<y_{j-1}$, then
\begin{equation*}
\Phi(\rho,\eta)
= \lam_1\otimes O(\rho) + \lam_2\otimes\bigl((m_{j},n_{j}) + 
O(\rho^2)\bigr).
\end{equation*}
This is indicated by the following diagram
\begin{diagram}[silent,width=2.2em,abut]
\null&\isot{(0,1)}{S_{k+2}}&\dotify{\mathop\bullet\limits_{y_{k+1}}^{p}}
&\isot{(p,q)}{S_{k+1}}&\dotify{\mathop\bullet\limits_{y_k}^{1}}
&\rLine&{\,\cdots\,}
&\rLine
&\dotify{\mathop\bullet\limits_{y_3}^{1}}
&\isot{\;\;(m_3,n_3)\;\;}{S_3}&\dotify{\mathop\bullet\limits_{y_2}^{1}}
&\isot{(m_2,n_2)}{S_2}&\dotify{\mathop\bullet\limits_{y_1}^{1}}
&\isot{(1,0)}{S_{1}}&\dotify{\mathop\bullet\limits_{y_{0}}^{1}}
&\isot{(0,-1)}{S_0}&\null
\end{diagram}
where the labelling conventions are as in \S\ref{sec:toric}.  In
particular, all points save $y_{k+1}$ are smooth points of the toric
manifold, while $y_{k+1}$ is an orbifold point.

It follows from this discussion and the results summarized from
Joyce~\cite{Joy:esd} in \S\ref{sec:cpt} that the conformal class defined by
$g_\Phi$ extends as a smooth orbifold metric to $\oX$.  Next one checks
that the `scalar-flat K\"ahler' conformal factor
\begin{equation} \label{1.19.2}
\Omega^2 = \frac{\rho|\det\Phi|}{\rho^2 + (\eta-y_{k+1})^2}
\end{equation}
has precisely the correct boundary behaviour so that $g_{SFK}$ of
\eqref{17.1.2} extends smoothly to $X$.

Finally we observe (as remarked by Joyce after the proof of~\cite[Theorem
3.3.1]{Joy:esd}) that the conformal factor~\eqref{1.19.2} ensures that the
metric is asymptotic near $y_{k+1}$ to the flat metric on $\R^4$ near
infinity. Since this is an orbifold singularity of $\oX$ in general, the
metric is ALE.  (When there is no singularity, we recover Joyce's half
conformally flat metrics on $k\CP2$.) The complex structure preserves the
special orbits, so it is associated to a realisation of $X$ as a resolution
of $\C^2/\Gamma$. This proves Theorem~\ref{thma}.

\vspace{5pt}

This result effectively gives a $k-1$ dimensional family of ALE scalar-flat
K\"ahler metrics, because an equivalent metric is obtained under a
projective transformation of the points $y_0,y_1,\ldots y_{k+1}$.  We could
equally well have assumed that the points are increasing (as does Joyce),
but to assume them decreasing is more natural in the next section.

\subsection{Self-dual Einstein metrics}\label{cansde}
We shall prove Theorem~\ref{thmb} using the results of~\cite{CaPe:emt}
summarized in Theorem~\ref{thm:cp}. We shall show that if all $e_j\geq 3$,
then the Joyce matrix $\Phi$ admits a potential, for a particular choice of
the $y_j$.  For any set of real numbers $w_j$, the eigenfunction
\begin{equation} \label{28.1.15}
F = \sum_{j=0}^{k+1} w_j F(\cdot ;y_j)
\end{equation}
is a potential for the Joyce matrix
\begin{equation}\label{28.1.16}
\Phi = \frac12\sum_{j=0}^{k+1} \phi(\cdot; y_j)\otimes (w_j,y_jw_j)
\end{equation}
and this has the form \eqref{28.1.2} if
\begin{equation*}
(m_{j}-m_{j+1},n_{j}-n_{j+1}) = (w_j, y_jw_j).
\end{equation*}
Reading this the other way, we see that the SFK metric on $X$ 
has (locally) a SDE representative $g_F$ in its conformal class, with 
potential given by \eqref{28.1.15} if and only if
\begin{equation*}
w_j = m_{j}-m_{j+1},
\end{equation*}
and the sequence
$y_0,\ldots y_{k+1}$ is given by
\begin{equation*}
y_j = \frac{n_{j+1} - n_{j}}{m_{j+1}-m_{j}}.
\end{equation*}
We require that this sequence is monotonic. Using the
definition~\eqref{2.13.3} of $e_j$ we obtain
\begin{equation*}
y_{j} - y_{j-1} = \frac{2-e_j}{(m_{j+1}-m_j)(m_{j}-m_{j-1})}
\end{equation*}
and if all $e_j\geq 2$, the sequence $m_j$ is strictly increasing.  It
follows that if all $e_j\geq 3$, then the sequence $y_j$ is strictly
decreasing.

\vspace{5pt}

The SDE representative is defined only where $F\neq 0$, so the next step is
to analyse the zero-set of $F$. Since
\begin{equation} \label{6.6.2}
f(\rho,\eta) := \sqrt{\rho}F(\rho,\eta)
\end{equation}
is continuous up to the boundary we may define
\begin{equation*}
Z = \{f =0\}, \quad
D_+ = \{f>0\}\quad\text{and}\quad D_- = \{f<0\}
\end{equation*} 
as subsets of the conformal 
compactification $\oH^2:=\{\rho\geq 0,\eta\in\R\}\cup\infty$ of $\cH^2$.

Notice that at any interior point of $Z$,
\begin{equation} \label{5.6.2}
|d F|^2 = - \det\Phi >0
\end{equation}
so that $F$ is a defining function and $Z\cap \cH^2$ is a closed
submanifold of $\cH^2$.  We claim now that $Z$ is a simple arc joining a
certain point $(0,\eta_0)$ to $\infty$ and intersecting no other point on
the boundary of $\cH^2$.  First we rule out the possibility that $Z$
contains a closed $C$ curve in the interior of $\cH^2$.  Because $F$ is a
defining function, such a curve is a smooth submanifold, and hence is {\em
simple}. By the Jordan curve theorem, $C$ has an interior $U$, say and
because $F$ is continuous on the compact set $U\cup C$, and is zero on the
boundary, it is either identically zero or achieves a positive maximum or a
negative minimum in $U$.  But as a solution of $\Delta F = (3/4)F$, $F$
satisfies a strong maximum principle which rules out the last two
possibilities. Hence $F=0$ in $U$, contradicting (for example)
\eqref{5.6.2}. We conclude that $Z$ consists of one or more disjoint simple
arcs joining boundary points of $\cH^2$.

In order to determine the end-points of $Z$, we must calculate the zeros of
$f(0,\eta)$.  It can be checked (see also Proposition~\ref{prop:bv}) that
if $y_{j}< \eta < y_{j-1}$, then
\begin{equation} 
f_0(\eta) := f(0,\eta)= m_{j}\eta-n_{j}.
\end{equation}
Since $m_{s}\geq 0$, $f_0(\eta)$ is non-decreasing in $(y_{s+1},y_s)$ and
since $f_0(\eta)$ is continuous, it follows that $f_0(\eta)$ is a
non-decreasing function of $\eta\in \R$.  Since $f_0(\pm\infty) = \pm 1$
and $f_0$ is strictly increasing in $[y_{k+1},y_0]$, $f_0(\eta)$ has a
unique zero at $\eta_0$, say (Figure~\ref{fcan}).
\ifaddpics
\begin{figure}[ht]
\begin{center}
\includegraphics[width=.7\textwidth]{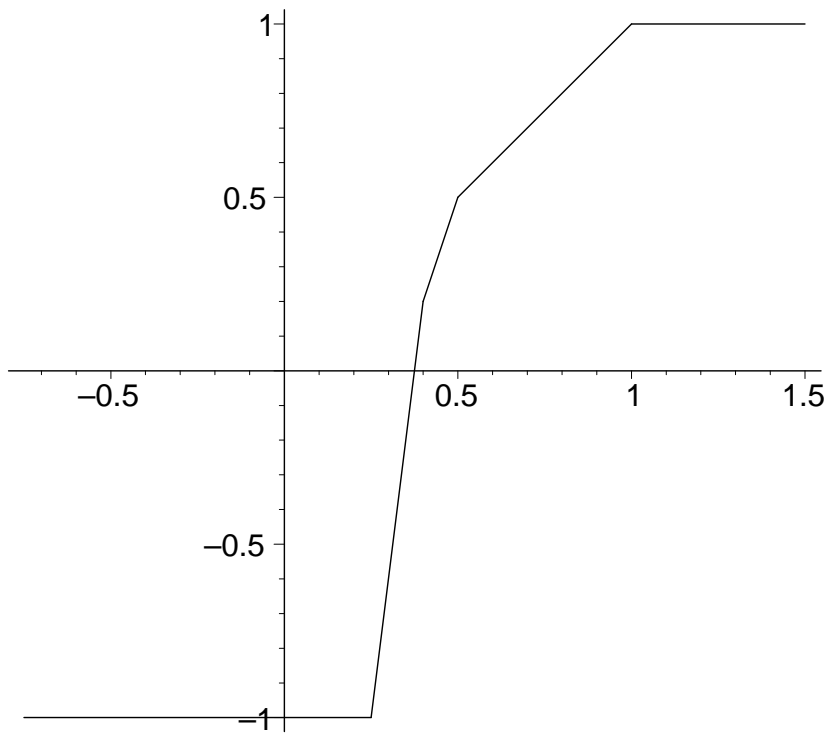}
\caption{The graph of $f_0$ when $q/p = 3/8$. The zero is at $3/8$.}
\label{fcan}
\end{center}
\end{figure}
\fi

In fact, since
\begin{equation*}
y_{k+1} = \frac{n_{k+2}-n_{k+1}}{m_{k+2}-m_{k+1}} = \frac{q-1}{p}
 < \frac{q}{p} <\frac{q-n_k}{p-m_k} =
\frac{n_{k+1}-n_{k}}{m_{k+1}-m_{k}} = y_k 
\end{equation*}
(note that $n_k/m_k < q/p$), it follows that $\eta_0= q/p$ and this point
lies in the interval $(y_{k+1},y_k)$.  Hence the boundary of $D_+$ is equal
to $(q/p,\infty)$.

Since $f$ is smooth up to the boundary near $(0,\eta_0)$, $Z$ must consist
of a single, simple arc joining this point to $\infty$.  Further, $\del_\rho
f(0,\eta_0)=0$ so $Z$ cuts the $\eta$-axis orthogonally. To complete the
proof, we simply define $X_\pm$ to be the set of orbits of $X$ lying over
$D_\pm$ and $Y$ to be the set of orbits lying over $Z$. Then the pull-back
of $f$ to $X$ is a defining function for $Y$ and the restriction $g_+$ of
$g_F$ to $X_+$ is the desired complete SDE metric.

The metric $g_F$ extends to the special orbits because $g_{SFK}$ has this
property and the ratio of the SDE and SFK conformal factors is
\begin{equation*}
\frac{\rho F^2}{\rho^2+(\eta -(q-1)/p)^2}
\end{equation*}
which has a continuous positive limit as $\rho\to 0$ if $\eta> q/p$. This
proves Theorem~\ref{thmb}.

\subsubsection{Remark} As mentioned in the introduction, the
restriction $g_-$ of $g_F$ to $X_-$ is also of interest.  It is not
difficult to check that $g_-$ extends to a smooth orbifold metric on $X_-$,
and passing to the universal cover we obtain a smooth AH $\Gamma$-invariant
SDE metric on the $4$-ball.

This picture was well known for the special case when $X$ reduces to
$\cO(-p)$. Let us now consider these Pedersen--LeBrun metrics in more
detail from this point of view.

\subsubsection{Example}
Here $k=1$ and there is only one parameter $p=e_1>0$.
We have
\begin{gather*}
(m_0,n_0) = (0,-1),\quad
(m_1,n_1) = (1,0),\quad (m_2,n_2) = (p,1),\quad (m_3,n_3)= (0,1)\\
\tag*{and}
y_0=1,\quad  y_1= 1/(p-1),\quad y_2=0.
\end{gather*}
We see at once that we must take $p\neq 2$ if the $y_j$ are to be distinct.
This fits with the fact that when $p=2$ the scalar-flat K\"ahler metric is
hyperk\"ahler.

If $p>2$, then we get examples of the type discussed here, with
\begin{equation*}
F = -F(\cdot;1)  - (p-1)F(\cdot;1/(p-1)) + pF(\cdot;0).
\end{equation*}

The case $p=1$ is slightly awkward from this point of view, but by applying
an orientation-preserving projective transformation, we see that we should
regard the points $1,\infty,0$ as increasing: indeed they are increasing
after a cyclic permutation.  The Einstein metric in this case has positive
scalar curvature: it is the Fubini--Study metric on $\CP2$ defined by
\begin{equation*}
F = F(\cdot;1) + F(\cdot; \infty) + F(\cdot;0)
\end{equation*}
where $F(\rho,\eta;\infty):= 1/\sqrt\rho$.

The results of this section show that exactly the same phenomenon occurs in
general, except that there are no further examples with positive scalar
curvature, and in the negative scalar curvature case we need $e_j>2$ for
all $j$.

%
% Proof of Theorem~\ref{thmc}
%
\section{Infinite-dimensional families of complete SDE metrics}
\label{sec:inni}

In this section, we shift viewpoint so that the `labelling' of the edges of
$M/T^2$ is regarded as a $\V$-valued step-function
\begin{equation*}
(m(y),n(y)) = (m_j,n_j)\quad \mbox{for}\quad y_j< y<y_{j-1},
\end{equation*}
where we regard $y_{k+2}= -\infty$, $y_{-1}= +\infty$.  Then the derivative
$(m',n')$ of $(m,n)$ with respect to $y$ is defined in the sense of
distributions and since
\begin{equation*}
(m',n') = \sum_{j=0}^{k+1} \delta(y-y_j)(m_j-m_{j+1},n_j-n_{j+1}),
\end{equation*}
we see that the formula~\eqref{28.1.2} for the matrix $\Phi$ used in the
last section can be written
\begin{equation*}
\Phi(\rho,\eta) = \frac12\int
\phi(\rho,\eta;y)\otimes(m',n')\,dy.
\end{equation*}
Note that we shall use `classical' notation for distributions (\ie, $\int
u(y)\phi(y)\,dy$ in place of $\ip{u,\phi}$ for the pairing of a
distribution $u$ with a test-function $\phi$) throughout this section.

The theme of this section is the replacement of the step-function $(m,n)$
by more general distributions.  This point of view gives a rather efficient
proofs of the results outlined in the previous section, in a much more
general setting.

\subsection{Smeared solutions of the Joyce equation and their
potentials}\label{sec:smear}

Let $v$ be a $\V$-valued distribution on $\R$ with compact support.
Then
\begin{equation} \label{3.1.2}
\Phi(\rho,\eta) = \frac 12 \int \phi(\rho,\eta;y)\otimes v(y)\,dy
\end{equation}
defines a smooth Joyce matrix in $\cH^2$, since $\phi(\rho,\eta;y)$ is
$\Cinf$ in all variables for $\rho>0$.  Similarly, if $w$ is any
real valued distribution on $\R$ with compact support, we may define
\begin{equation} \label{4.1.2}
F(\rho,\eta) = \int F(\rho,\eta;y)w(y)\,dy
\end{equation}
and this will be a smooth eigenfunction of the hyperbolic laplacian. These
are the `smeared' solutions of the title.  Moreover, exactly as in
\S\ref{cansde}, \eqref{4.1.2} is a potential for \eqref{3.1.2} if we set
\begin{equation} \label{5.1.2}
v(y) = (w(y), y w(y)).
\end{equation}

Motivated by the introduction to this section, we wish to write
$v=(\mu',\nu')$ where $\mu$ and $\nu$ are distributions on $\R$.  Recall
that a distribution $u$ satisfies $u'=0$ (in the sense of distributions) on
an open set if and only if $u$ is locally constant on this open set.
Therefore $v$ has compact support if and only if
\begin{equation}\label{2.1.2}
\mu\text{ and }\nu\text{ are {\em locally constant near infinity}}.
\end{equation}
We shall assume additionally that $\mu$ and $\nu$ are {\em odd at
infinity}, \ie,
\begin{equation} \label{55.1.2}
\mu(\infty) = - \mu(-\infty),\qquad\nu(\infty) = - \nu(-\infty).
\end{equation}

Suppose now that~\eqref{5.1.2} also holds, \ie,
\begin{equation}
\mu'(y)=w(y),\qquad \nu'(y) = yw(y),
\end{equation}
and set
\begin{equation} \label{7.17.3}
f_0(y)=y\mu(y)-\nu(y).
\end{equation}
Then $f_0''(y)=w(y)$ and 
\begin{equation}\label{5.1.3}
\mu(y)=f_0'(y),\qquad \nu(y) = yf_0'(y)-f_0(y).
\end{equation}
Conversely, if $f_0$ is a distribution on $\R$ such that $f_0''$ has
compact support, then \eqref{7.17.3} holds with $\mu$ and $\nu$ defined
by~\eqref{5.1.3}; $f_0$ is {\em affine near infinity}, and we also require,
in accordance with~\eqref{55.1.2}, that the coefficients $\mu,\nu$, are
not just locally constant, but odd at infinity.

These assumptions allow us to integrate by parts once in \eqref{3.1.2} and
twice in \eqref{4.1.2} to get the formulae
\begin{equation} \label{5.17.3}
\Phi(\rho,\eta) = 
\frac{\rho}{2}\int
\frac{(y-\eta)\lam_1+\rho\lam_2}{\bigl(\rho^2 + (\eta-y)^2\bigr){}^{3/2}}
\otimes(\mu(y),\nu(y))\,dy
\end{equation}
and
\begin{equation} \label{6.17.3}
F(\rho,\eta) = \frac{\rho^{3/2}}{2}\int\frac{f_0(y)\,dy}
{\bigl(\rho^2 + (\eta-y)^2\bigr){}^{3/2}}.
\end{equation}
These are Poisson formulae in the following sense.

\begin{prop} \label{prop:bv}
If $\Phi$ and $F$ are given by \eqref{5.17.3} and
\eqref{6.17.3} then we have \textup(in the sense of distributions\textup)
\begin{equation}\label{5.19.3}
\Phi(\rho,\eta) \to \lam_2\otimes(\mu(\eta),\nu(\eta))
\mbox{ as }\rho\to 0
\end{equation}
and
\begin{equation} \label{6.19.3}
\sqrt{\rho}F(\rho,\eta) \to f_0(\eta)\mbox{ as }\rho\to 0.
\end{equation}
Moreover, if $(\mu,\nu)$ is constant and $f_0$ is affine
in a neighbourhood of $\eta = a$ then
\begin{equation}\label{7.19.3} 
\Phi(\rho,\eta) = \lam_1\otimes O(\rho)+ \lam_2\otimes
\bigl((\mu(\eta),\nu(\eta))+ O(\rho^2)\bigr)
\end{equation}
and
\begin{equation}\label{1.31.5}
\sqrt\rho F(\rho,\eta) = f_0(\eta) + O(\rho^2)
\end{equation}
near $\eta=a$.
\end{prop}

\subsubsection{Remark} From now on, when we say that an eigenfunction
$F$ has a boundary value $f_0$, we shall always mean \eqref{6.19.3}.

\subsubsection{Proof} We claim that 
\begin{equation*}
\lim_{\rho\to 0}\frac{\rho^2}{2(\rho^2 +y^2)^{3/2}} = \delta(y)
\end{equation*}
and that
\begin{equation*}
\lim_{\rho\to 0}\frac{y}{2(\rho^2 +y^2)^{3/2}} = 0.
\end{equation*}
These two limits yield \eqref{5.19.3} and \eqref{6.19.3} (by
translation of $y$).  To establish the first claim, test against a function
$\phi$ and make the change of variables $y =\rho z$ to get
\begin{equation*}
\int \frac{\rho^2}{2(\rho^2 +y^2)^{3/2}}\phi(y)\,dy = 
\int \frac{\phi(\rho z)\,dz}{2(1 +z^2)^{3/2}} \to \phi(0)
\end{equation*}
as $\rho \to 0$,  since $\int (1+z^2)^{-3/2}\,dz = 2$. Similarly, for the
second claim,
\begin{equation*}
\int \frac{\rho y }{2(\rho^2 +y^2)^{3/2}}\phi(y)\,dy = 
\int \frac{z \phi(\rho z)\,dz}{2(1 +z^2)^{3/2}} \to 0
\end{equation*}
as $\rho \to 0$. 

If now $(\mu,\nu)$ is constant near $\eta=a$ we can expand in powers of
$\rho^2$ the denominator $\bigl(\rho^2+(a-y)^2\bigr){}^{-1/2}$ in
\eqref{3.1.2}.  Such an expansion is valid provided $\rho$ is smaller than
the distance from $a$ to the support of $(\mu',\nu')$. This gives
\eqref{7.19.3}. A similar argument, starting from \eqref{4.1.2}, and
expanding $F(\rho,\eta;y)$ in powers of $\rho$, yields \eqref{1.31.5}.

\subsubsection{Remark} $\Phi$ and $F$ have complete asymptotic
expansions as $\rho\to 0$ whatever the behaviour of $(\mu,\nu)$ or $f_0$,
but in general, the next term in the expansion of $\sqrt\rho F$ is
$O(\rho^2\log\rho)$.

\subsubsection{Remark} It is important to ask what happens at $\infty$.
It would be quite straightforward to analyse this separately, but it is more
geometric to observe that all quantities are restrictions of globally
defined objects on the boundary of $\cH^2$. For example, if we set
\begin{equation*}
\trho = \frac{\rho}{\rho^2+\eta^2},\quad
\teta = -\frac{\eta}{\rho^2+\eta^2},\quad
\ty = -\frac{1}{y}
\end{equation*}
then
\begin{equation*}
\frac{\rho^2 + (\eta - y)^2}{\rho}(dy)^{-1} =
\frac{\trho^2 + (\teta - \ty)^2}{\trho}(d\ty)^{-1}.
\end{equation*}
These formulae are standard in analysis on hyperbolic space and imply that
the analysis near $\rho= \infty$ can be replaced by the analysis near
$\trho=0$ which has already been done. In particular, we see that
$F(\rho,\eta;y)$ has an invariant interpretation over
$\overline\cH{}^2\times \RP1$ as a section (with singularities at a certain
subset of the boundary) of $\cO(1)\tens L$, where $L$ is the M\"obius line
bundle over $\RP1$.  Similarly, $\phi(\rho,\eta;y)$ is a section of
$\cW^*\otimes L$. It is the trivialization of $L$ over
$\RP{1}\setdif\{\infty\}$ that corresponds to the `odd at infinity'
condition that is so prominent in this section. This global interpretation
allows us to relax the assumption that our distributions have compact
support on $\R$, so long as they are distributions at infinity, and we have
already taken advantage of this in our Poisson formulae. For more details,
see the appendix.

\vspace{5pt}

We now come to a useful sufficient condition for $\det\Phi\neq 0$.

\begin{prop} \label{prop:neg}
If $\Phi$ is given by \eqref{5.17.3} and
\begin{equation} \label{10.17.3}
\mu(y)\nu(z) -\mu(z)\nu(y) \leq 0\mbox{ for }y\leq z
\end{equation}
with strict inequality for some $y<z$, then $\det\Phi <0$ in $\cH^2$.
\end{prop}

Note that \eqref{10.17.3} makes good sense in terms of the
tensor product of distributions even when $\mu$ and $\nu$ are not
continuous functions.

\subsubsection{Proof} Recall that with our orientation conventions, if
$\Phi = \lam_1\otimes v_1 + \lam_2\otimes v_2$
then $\det\Phi = -  \ip{v_1,v_2}$.
Hence,
\begin{equation} \label{2.19.2}
\det\Phi(\rho,\eta) = -\frac{\rho^3}{8}\iint
\frac{(y-z)(\mu(y)\nu(z)-\mu(z)\nu(y))}
{\bigl(\rho^2+(\eta-y)^2\bigr){}^{3/2}\bigl(\rho^2+(\eta-z)^2\bigr){}^{3/2}}
\,dy\,dz.
\end{equation}
The result follows at once.

\subsubsection{Remark} If $\mu$ and $\nu$ are piecewise continuous
functions with $\mu\geq 0$, then \eqref{10.17.3} is equivalent to
\begin{equation} \label{3.19.2}
\frac{\nu(y)}{\mu(y)}\geq \frac{\nu(z)}{\mu(z)}
\quad \mbox{for}\quad y\leq z
\end{equation}
so that $\nu(y)/\mu(y)$ is a  non-increasing function of $y$.  This property
is enjoyed by the step-function $(m,n)$ and is used in Joyce's
proof of the non-vanishing of $\det\Phi$.

\subsubsection{Remark} We have now associated to any pair of
distributions $(\mu,\nu)$ on $\R$ that satisfy \eqref{2.1.2},
\eqref{55.1.2} and \eqref{10.17.3} a $\V$-invariant conformally half-flat
metric $g$ on $\cH^2\times \V$.

\subsection{Infinite-dimensional families of SDE metrics}

The $\V$-valued step function $(m,n)$ associated to a Hirzebruch--Jung
resolution $X$ with $c_1(X)<0$ gave rise to an eigenfunction $F$ with
boundary data
$$
f^{can}_0(\eta) = m_j\eta -n_j\quad \mbox{for}\quad y_j\leq \eta \leq
y_{j-1}
$$
where $y_{k+2}=-\infty$ and $y_{-1}= +\infty$ as before.  In particular
$f^{can}_0(\eta) = \pm 1$ for $\pm \eta$ sufficiently large and positive.
The graph of $f^{can}_0$ is continuous, strictly increasing in
$[y_{k+1},y_0]$ and has the property that it is {\em convex} where it is
negative and {\em concave} where it is positive. The unique zero of
$f^{can}_0$ is at $\eta = q/p$. (Recall the example of Figure~\ref{fcan}.)

In order to generate an infinite-dimensional family of complete SDE metrics
on $X$, we shall change $f_0^{can}$ in the interval $(-\infty,q/p)$.  (If
we change it in $(q/p,\infty)$, the result may no longer extend smoothly to
the special orbits over the intervals $[y_j,y_{j-1}]$.) For our first
statement, we consider modifications $f_0$ that are not too rough.

\begin{thm}\label{thm:inf1}
Let $f_0(\eta)$ be a continuous function equal to $f_0^{can}(\eta)$
for $\eta\geq q/p-\delta$, for some $\delta>0$.  Suppose further that
$f_0(\eta)=-1$ for all $\eta\leq a<q/p$ and that $f_0$ is piecewise
differentiable, strictly increasing and convex on $[a,q/p]$. Then the
eigenfunction $F$ with boundary value $f_0$ determines an
asymptotically hyperbolic SDE metric on the neighbourhood $E\subset
X_+\subset X$ of the exceptional divisor corresponding to the domain
$F>0$ in $\cH^2$.
\end{thm}

\subsubsection{Proof} Define distributions $\mu,\nu$ by~\eqref{5.1.3}.
Then $(\mu,\nu)$ determines a Joyce matrix $\Phi$ by \eqref{5.17.3} with
the correct boundary values for smooth extension to the special orbits
$S_1,S_2,\ldots S_k$. Let us check that $\det\Phi$ is never zero. Since
$\mu$ is defined and is positive a.e.\ in $[a,y_0]$, we can apply
Proposition~\ref{prop:neg} in the form of \eqref{3.19.2}. Indeed,
\begin{equation*}
\frac{\nu(y)}{\mu(y)} = y - \frac{f_0(y)}{f_0'(y)}
\qquad\text{and so}\qquad
\Bigl(\frac{\nu}{\mu}\Bigr)' =   \frac{f_0f_0''}{(f_0')^2} \leq 0
\end{equation*}
by the assumption that $f_0$ is concave where $f_0>0$ and $f_0$ is 
convex where $f_0<0$. Hence $\nu/\mu$ is a non-increasing function as 
required.  Thus the conformal class of $g_\Phi$
is defined over the whole of $T^2\times \cH^2$. We now proceed as in
\S\ref{cansde}.

By the maximum principle, $F$ has no interior positive maximum or 
negative minimum. Hence the zero-set $Z$ is a simple smooth arc 
joining the boundary at $\eta= q/p$ to $\infty$ and $Z$ decomposes
$\oH^2$ into pieces
$$
D_{+} =\{ F >0\},\;\; D_{-} =\{ F <0\}.
$$
We define $X_+$ to be the union of $T^2$-orbits over $D_+$ and $Y$ to be
the union of the $T^2$-orbits over $Z$. We note (as before) that $X_+$
contains $S_1\cup \cdots \cup S_k$ (since $f_0$ is positive to the right of
$\eta=q/p$). Then the metric $g_F$ defines an asymptotically hyperbolic SDE
metric on $X_+$, with conformal infinity $Y$.

\subsubsection{Remark} In the situation of Theorem~\ref{thm:inf1}
we can say more about the zero-set $Z$ of $F$. We claim that
it meets
each circular arc with end-points $(0,q/p \pm b)$ in precisely one
point, for any $b>0$.
To see this, note that these arcs are the orbits of the Killing
vector field
\begin{equation*}
K = \bigl((\eta-q/p)^2 -b^2\bigr)\del_\eta + 2(\eta-q/p)\rho\,\del_\rho.
\end{equation*}
Then $K\cdot F$ will be an eigenfunction of the laplacian, with eigenvalue
$3/4$ and so by the maximum principle, $K\cdot F<0$ in $\cH^2$ if this is
true near the boundary. Hence if we prove the latter, it will follow that
$F$ is strictly decreasing along these circular arcs, and since it starts
positive and ends negative, there must be a unique zero.  By
Proposition~\ref{prop:bv},
\begin{equation*}
\sqrt\rho F(\rho,\eta) \simeq \eta\mu(\eta)-\nu(\eta)\quad\mbox{if $\rho$
is small.}
\end{equation*}
Hence
\begin{align*}
\sqrt\rho\, K\cdot F &\simeq (\eta - q/p)^2\mu(\eta) -
(\eta-q/p)(\eta\mu(\eta)-\nu(\eta))\\
 &=\frac{1}{p}\bigl( -b^2\mu(\eta)  + (\eta-q/p)
(\mu(q/p)\nu(\eta) -\mu(\eta)\nu(q/p)\bigr) \leq 0
\end{align*}
by considering separately the cases $\eta <q/p$ and $\eta > q/p$ and 
using $\mu(q/p) = p$ and $\nu(q/p) = q$ as well as the
monotonicity property \eqref{10.17.3}.

\vspace{5pt} Theorem~\ref{thm:inf1} produces a family of SDE metrics
parameterized by piecewise differentiable functions on $(-\infty,q/p)$
satisfying only monotonicity and convexity assumptions. It is clear that
the space of such functions is (continuously) infinite-dimensional. Hence
Theorem~\ref{thm:inf1} implies Theorem~\ref{thmc}.  Notice that this result
is not a perturbation theorem: to the left of $q/p$, $f_0$ can be far from
$f_0^{can}$. A particularly interesting example is as follows.

\subsubsection{Example} Following an initial suggestion of Atiyah,
we note that the odd extension, $f_0^{odd}$, say, of $f_0^{can}$ to the
left of $\eta=q/p$ satisfies all the properties of
Theorem~\ref{thm:inf1}---see Figure~\ref{fodd}. More precisely, for
$\eta>0$, we put
$$
f_0^{odd}(q/p-\eta) = - f_0^{can}(q/p+\eta)
$$
so that
that $(\mu,\nu)$ satisfy
\begin{equation} \label{1.14.2}
\mu(q/p-\eta) = \mu(q/p+\eta),\quad \nu(q/p-\eta) =
2q\mu(q/p+\eta)/p-\nu(q/p+\eta).
\end{equation}
This extension enjoys the property that $F(\rho,q/p+\eta) =
-F(\rho,q/p-\eta)$ so that $Z= \{\eta = q/p\}$ and $D_+= \{\eta>q/p\}$.  In
this case we can compute a representative metric $h$ for the conformal
infinity by substituting $\eta=q/p$ into the formula for $F^2 g_F$. Since
$F$ vanishes along $\eta=q/p$, so does its tangential derivative
$F_\rho$---hence the conformal infinity is determined by the function
$F_\eta(\rho,q/p)$, which is a sum of the terms
\begin{equation*}
F_\eta(\rho,q/p;y_j)-F_\eta(\rho,q/p;2q/p-y_j)
= \frac{2(q/p-y_j)}{\sqrt\rho\sqrt{\rho^2+(q/p-y_j)^2}}
\end{equation*}
over $j=0,\ldots k$ (where $y_j>q/p$). We now readily obtain
\begin{equation*}
h=\biggl(\sum_{j=0}^k \frac{2(y_j-q/p)\sqrt{\rho}}
{\sqrt{\rho^2+(y_j-q/p)^2}}\biggr)^2
\frac{d\rho^2}{\rho^2}
+\rho \,d\phi^2 -\frac1\rho (d\psi-q/p \, d\phi)^2.
\end{equation*}

By definition, the odd extension has an obvious symmetry about $\eta=q/p$
and so $-f_0^{odd}$ should define a toric self-dual Einstein metric on a
manifold $X_-$ diffeomorphic to $X_+$. However, to do this, a different
integral lattice is needed to define the torus---otherwise $X_-$ will not be
smooth. More precisely, the lift of the symmetry $y\mapsto 2q/p - y$ to a
matrix of determinant $-1$ in $\GL_2(\R)$ is
\begin{equation*}
\begin{pmatrix} 1 & 0 \\ 2q/p & -1
\end{pmatrix} =
\begin{pmatrix} p & 0 \\ q & -1
\end{pmatrix}
\begin{pmatrix} p & 0 \\ q & 1
\end{pmatrix}^{-1}
\end{equation*}
and the new lattice is the image of $\Z^2$ under this transformation.  The
two lattices have a common sublattice of index $p$, \ie, away from the
fixed points there is a common covering space. We deduce that for each
choice of orientation, $h$ defines a conformal structure on $S^3$ with a
$\Z_p$-quotient bounding a self-dual Einstein metric; however, the two
$\Z_p$ actions are \emph{not} the same!

\ifaddpics
\begin{figure}[ht]
\begin{center}
\includegraphics[width=.7\textwidth]{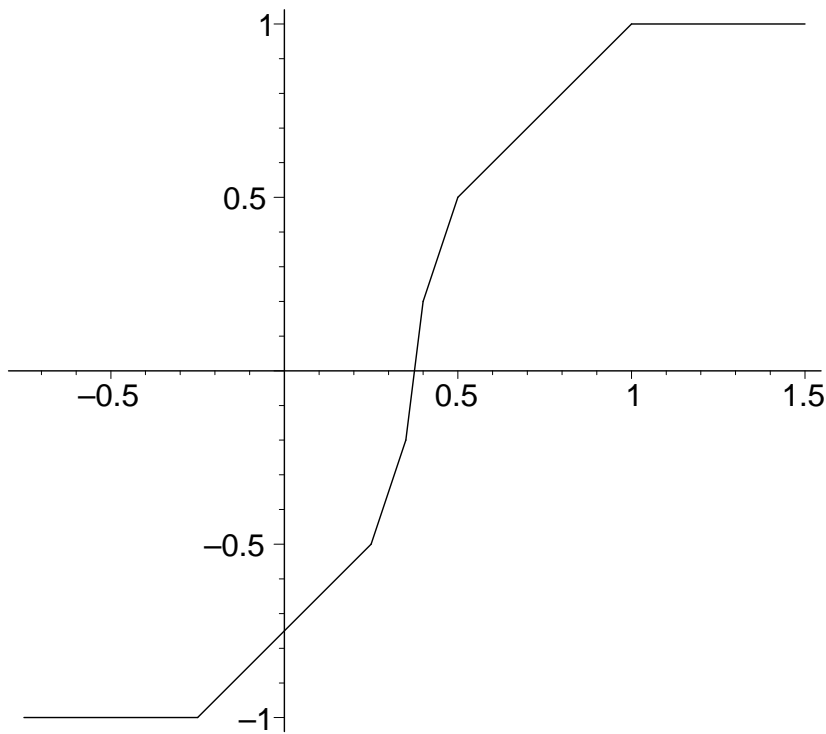}
\caption{The graph of $f^{odd}_0$ when $q/p = 3/8$.}
\label{fodd}
\end{center}
\end{figure}
\fi

\subsection{Further results}

The metrics constructed by Theorem~\ref{thm:inf1} all have the property
that the underlying self-dual conformal structures extend `a long way' into
$X_-$. More precisely they extend to the complement of the special orbits
(a set of codimension $2$) in $X_-$.  The reason for this is that we
arranged $\det\Phi\neq 0$ over the whole of $\cH^2$.

One would expect, however, that there should exist asymptotically
hyperbolic SDE metrics on $X_+$, with conformal infinity at $Y$, but with
the property that the self-dual conformal structure does not extend so far,
or even at all into $X_-$!

\subsubsection{Perturbations \textup I} Let $f_0$ be a function 
satisfying the properties of Theorem~\ref{thm:inf1} and let $u$ be 
any distribution with compact support contained in $(-\infty,b)$, 
say, where $b<q/p$. Consider, for real $t$, 
$$
f_0^{t} = f_0 + tu.
$$
For small $t$, the eigenfunction $F^t$ with boundary value $f_0^t$ will be
close to $F$, in the $\Cinf$-topology, on any set of the form
$$
\oH^2\setdif\{(\rho,\eta):0\leq\rho\leq\rho_0\mbox{ and }\eta\leq b\}.
$$
In particular, for all sufficiently small $t$, the zero-set of $F^t$ 
will be very close to that of $F$, and $\det\Phi^t\not=0$ where 
$F^t\geq0$. Hence $F^t$ yields a deformation of $g_F$ as an 
asymptotically hyperbolic SDE metric.

\vspace{5pt} Note that there is now no reason why $\det\Phi^t$ should be
nonzero on {\em all} of $\cH^2$. The $3$-pole examples given
in~\cite{CaPe:emt} illustrate this point.

\subsubsection{Perturbations \textup{II}}

A more general perturbation can be obtained from the following 
analytical result \cite{RM:pc}.

\begin{thm} Let $Z$ be a simple closed curve dividing $\cH^2$ into two
connected components $D_\pm$. Denote by $\del\cH^2_\pm$ the two 
components of the boundary of $\cH^2$, so that
\begin{equation*}
\del D_\pm = Z\cup\del\cH^2_\pm.
\end{equation*}
Then there exists a unique eigenfunction $F$ \textup(with eigenvalue
$3/4$\textup) on $D_+$ with the property that $F=0$ on $Z$ and $F$ has
prescribed boundary value on $\del \cH^2_+$\textup:
\begin{equation*}
f_0(\eta) = \lim_{\rho\to 0} \sqrt{\rho}F(\rho,\eta).
\end{equation*}
\end{thm}

To use this result, let $Z_0$ be the zero-set of an eigenfunction 
$F_0$ with boundary value $f_0$ satisfying the conditions of 
Theorem~\ref{thm:inf1}. Let $Z$ be a small perturbation of $Z_0$ (with 
end-points fixed) and let $F$ be the eigenfunction, with the same 
boundary values as $F_0$ for $\eta\geq q/p$ and vanishing on $Z$. If 
$Z$ is sufficiently close to $Z_0$ then $F$ will be close to $F_0$ 
in $D_+$, and so $F^2-4|dF|^2$ will not vanish on $D_+$. Accordingly 
the metric $g_F$ will be an asymptotically hyperbolic SDE metric with 
conformal infinity on $Z$.

Notice that if $F$ extends beyond its zero-set $Z$, then the latter will be
a real-analytic curve, since $F$ is real-analytic in the interior of its
domain of definition. Hence if we choose $Z$ to be merely smooth, $F$ will
not extend beyond $Z$ and the self-dual conformal structure of $g_F$ cannot
be extended through $Z$. These examples fill in non-analytic conformal
structures on lens spaces, cf.~\cite{Biq:mab}.

\subsubsection{Proof of Theorem~\ref{thmch}}

Consider the extension $f_0$ of $f_0^{can}$ by zero to the left of
$\eta=q/p$. This does not satisfy all of the conditions of
Theorem~\ref{thm:inf1}: in particular it is not odd at infinity. In
invariant terms (over $\RP1$), this means that $f_0$, though
continuous, is not smooth at infinity, but has a corner.  One way to handle
this is to change coordinates so that a smooth point of $f_0$ is at
infinity. However, it is straightforward to work directly with the given
coordinates, and this is what we shall do.

It is immediate that $\det\Phi<0$ on $\cH^2$. Furthermore $F$ is
positive on $\cH^2$ by the Poisson formula, so the Einstein metric is
smoothly defined on {\it all} of $\cH^2\times T^2$. Since $f_0$ agrees
with $f_0^{can}$ on $(q/p,\infty)$, the metric extends smoothly to the
special orbits over this interval. It remains to consider the
behaviour of the metric as $(\rho,\eta)$ approaches the boundary
segment $[-\infty,q/p]$.

We claim that this is a complete end of the self-dual Einstein metric and
that it is ACH.  To see this, we compare $f_0$ to the function
\begin{equation*}
f_0^{CH}=\begin{cases} 0 & \text{if }\eta<q/p\\
p\eta-q & \text{if }q/p<\eta<(q+1)/p\\ 1 & \text{if }(q+1)/p<\eta.
\end{cases}
\end{equation*}
We observe that $f_0$ and $f_0^{}$ agree, except on
$[y_k,y_0]\subset (q/p,\infty)$. On the other hand $f_0^{CH}=
\frac12(|p\eta-q|-|p\eta-q-1|+1)$ and so it generates a $3$-pole solution
in the sense of~\cite{CaPe:emt}.

In fact after change of $(\rho,\eta)$ coordinates we find that
\begin{equation*}
F^{CH}(\rho,\eta) = \sqrt{\frac p2}\biggl(
-\frac{1}{\sqrt\rho}
+\frac{\sqrt{\rho^2+(\eta+1)^2}}{2\sqrt\rho}
+\frac{\sqrt{\rho^2+(\eta-1)^2}}{2\sqrt\rho}\biggr),
\end{equation*}
which is a hyperbolic eigenfunction generating the Bergman
metric~\cite{CaPe:emt}.  (Explicitly, with $\rho=2\coth t\,{\rm
csch}\,t\,\sin\theta$, $\eta=(2\coth^2 t-1)\cos\theta$, we have
\begin{equation*}
g_{CH}= 2dt^2 + \frac12 \sinh^2 t\, ( d\theta^2+ (2/p)\sin^2\theta\,
d\phi^2) + \frac1{4p} \sinh^2 2t\, (d\psi+\cos\theta\,d\phi)^2
\end{equation*}
defined on a $\Z_p$ orbifold quotient of $\C\cH^2$. Rescaling $g_{CH}$ by
$1/2$, and $\phi$ and $\psi$ by $\sqrt{p/2}$ gives the form of the Bergman
metric given in~\cite{Rol:reh}.)

The approximation of $f_0$ by $f_0^{CH}$ is enough to
ensure completeness. In particular, one can easily check directly that
$\eta=q/p$ and $\eta=\pm\infty$ are `at infinity' on the special orbits
over $(q/p,y_k)$ and $(y_0,\infty)$.

We now compute the CR infinity, by taking $\eta\in (-\infty,q/p)$ and the
limit $\rho\to 0$. Since $f_0=0$, and hence $(\mu,\nu)=(0,0)$, on
$(-\infty,q/p)$, we can make an asymptotic expansion for $\rho$ smaller
than $q/p-\eta$ to obtain
\begin{align*}
\sqrt\rho F(\rho,\eta)&=\tfrac12 \rho^2 f_1(\eta)+O(\rho^4)\\
v_1&= \rho(f_1(\eta),\eta\,f_1(\eta))+O(\rho^3)\\
v_2&= \tfrac12\rho^2\bigl(f_1'(\eta),f_1(\eta)+\eta\,f_1'(\eta)\bigr)
+O(\rho^4)
\end{align*}
where
\begin{equation}
f_1(\eta)=\int \frac{f_0''(y)}{|\eta-y|} dy.
\end{equation}
(Recall that $f_0''$ is a sum of delta distributions.)

We let $\theta$ be the $1$-form $\rho^4\ip{v_1,\cdot}^2/F^2\ip{v_1,v_2}$
and note that $\lim_{\rho\to 0}\theta =
2(d\psi+\eta\,d\phi)^2/f_1(\eta)^2$. Now $h=\rho^2 g-\rho^{-2}\theta^2$
is degenerate, and we compute that
\begin{equation*}
\lim_{\rho\to 0} h =\frac{4\,({d\eta}^2 + {d\rho}^2) \,{f_1(\eta)}^4 + 
    {( f_1'(\eta)d\psi + (f_1(\eta) + \eta\,f_1'(\eta)) d\phi)}^2}
{2\,{\rho}^2{f_1(\eta)}^4}.
\end{equation*}
Thus, after rescaling by $f_1(\eta)^2$, the pullback of $(\theta,h)$ to
$\rho=0$ (restricting $h$ to $\ker\theta$) gives the contact metric
structure
\begin{equation}\label{CRstr}
2(d\psi+\eta\,d\phi),\qquad 2f_1(\eta)d\eta^2+\frac{d\phi^2}{2f_1(\eta)}
\end{equation}
The Reeb field $\partial_\psi$ is CR and generates the foliation of the
lens space $N$ induced by the Hopf fibration of $S^3$, and so $N$ is normal
and quasiregular. The quotient metric is an $S^1$-invariant orbifold metric
on $S^2$, as one can easily check directly from~\eqref{CRstr}.

\section{Proof of Theorem~\ref{thmd}}

Recall from previous sections that for an eigenfunction $F$ defined 
in $\cH^2$, the metric $g_F$ extends to the special orbits over the 
interval $(a,b)$, with isotropy $G(m,n)$ if
$$
f(\rho,\eta) = f_0(\eta) +\frac{1}{2}\rho^2\,f_1(\eta)+\cdots
\quad\mbox{for $a<\eta<b$ and $\rho$ small},
$$
where
$$
f_0(\eta) = m\eta -n.
$$
Furthermore, $g_F$ extends smoothly through the fixed-points
corresponding to $a$ and $b$ if $f$ has similar expansions for $\eta$ to
the left of $a$ and to the right of $b$, and $f_0(\eta)$ is continuous at
$a$ and $b$.  As a warm-up for the discussion in the rest of this section,
let us sketch a proof of this.  If $f$ has an expansion as above, then to
leading order
$$
v_1 = (\rho f_1,\rho(\eta f_1-m)),\;\;v_2 = (m,n)
$$
and
$$
\ip{v_1,v_2} = \rho(n f_1 - m(\eta f_1-m)) = \rho(m^2-f_0f_1).
$$
To leading order, therefore,
$$
g_F = \frac{|m^2-f_0f_1|}{f_0^2}\left(
d\rho^2 +d\eta^2 + \frac{(nd\theta_1-md\theta_2)^2 + 
\rho^2(f_1d\theta_2 +(m-\eta f_1)d\theta_1)^2}{(m^2-f_0f_1)^2}\right).
$$
Pick integers $a$, $b$ with $an-mb=1$ and make the change of basis
$$
d\theta_1 = a\,d\phi + m \,d\psi,\quad d\theta_2 = b\, d\phi + n\, d\psi
$$
so that 
$$
n \, d\theta_1 - m \, d\theta_2 = d\phi,\quad
f_1 d\theta_2 +(m-\eta f_1) d\theta_1
= (m^2- f_0f_1)d\psi +(am - f_1(a\eta-b) )d\phi.
$$
Then
$$
g_F = \frac{|m^2-f_0f_1|}{f_0^2}\left(
d\rho^2 + \rho^2d\psi^2 + d\eta^2 + \mbox{terms in $d\phi^2$ and 
$d\phi \, d\psi$}\right).
$$
In particular the metric in the $(\rho,\psi)$-plane is a smooth 
rescaling of the flat metric in polar coordinates, and hence extends 
to $\rho=0$.

The `corners' are dealt with similarly, using the change of variables $\rho
= 2r_1r_2$, $\eta = r_1^2-r_2^2$.

\subsection{Outline of proof}\label{outline}
We know that every toric SDE metric arises from a suitable function $F$. We
know also that an embedded sphere consisting of special orbits must
correspond to a combinatorial diagram of the form
\begin{diagram}[silent,width=2.2em,abut]
\null
&
\isot{\;\;(m,1)\;\;\;\;}{S_{2}}&\dotify{\mathop\bullet
\limits_{y_1}^{1}}&
\isot{\;\;\;(1,0)\;\;\;}{S_1}&
\dotify{\mathop\bullet\limits_{y_0}^{1}}&
\isot{\;\;\;\;(0,-1)\;\;}{S_{0}}&
\null
\end{diagram}
where $m$ is a positive integer and we have used $SL(2,\Z)$ invariance to
fix the isotropies of $S_0$ an $S_1$. The embedded sphere here is
$S_1$---it has self-intersection $m$, and so we must show that $m\geq 3$ to
prove Theorem~\ref{thmd}.  The eigenfunction $F$ is defined in a one-sided
neighbourhood $U$ of this diagram and $U$ carries a metric $h$ of constant
curvature $-1$.  The boundary of $U$ will be divided into two pieces, $\del
U = C_1\cup C_2$ such that $C_1$ contains $S_1$ in its interior.  Suppose
next that we have already proved that if $g_F$ extends to the special
orbits, then it must do so in the `standard way'; that is, the boundary of
the diagram is part of the conformal infinity of $\oH^2$ and $f$ has the
above standard form near the boundary:
\begin{equation} \label{2.18.3}\begin{split}
f_0(\eta) &=
\begin{cases}
m\eta - 1 & \mbox{ if }\eta\leq y_1\\
\eta & \mbox{ if }y_1<\eta <y_0 \\
1 & \mbox{ if }y_0\leq\eta
\end{cases}\\
\mbox{and }f_0&\mbox{ is continuous at $y_0$ and $y_1$.}
\end{split}\end{equation}
For $m\neq 1$, the continuity of $f_0$ means that
\begin{equation*}
y_1 = 1/(m-1),\quad y_0=1.
\end{equation*}
The result now follows because if $m=2$ then $y_0=y_1$, which is not
allowed, and if $m<1$, then $y_1<0$ and $F$ changes sign near $S_1$, so the
corresponding SDE metric blows up there. The case $m=1$ is a little tricky
with this normalization, reflecting a poor choice of basis (or point at
infinity) in this case. However if we switch to
\begin{diagram}[silent,width=2.2em,abut]
\null&
\isot{\;\;(1,1)\;\;\;\;}{S_{2}}&\dotify{\mathop\bullet\limits_{y_1}^{1}}&
\isot{\;\;\;(2,1)\;\;\;}{S_1}&
\dotify{\mathop\bullet\limits_{y_0}^{1}}&
\isot{\;\;\;\;(1,0)\;\;}{S_{0}}&
\null
\end{diagram}
then continuity gives $y_1=0$, $y_0=1$ and
\begin{equation*}
f_0(\eta) =
\begin{cases}
 \eta -1 &  \mbox{ if }\eta\leq  0\\
2\eta-1 & \mbox{ if } 0<\eta < 1 \\
\eta &  \mbox{ if }1\leq \eta
\end{cases}
\end{equation*}
so that $f$ changes sign in the middle of $S_1$. Hence this cannot
happen for the same reason as before.  It is good to note that we do
not rule out the existence of the standard metrics on $S^4$ and $\C
P^2$ this way. The reason for this is that those metrics have positive
scalar curvature, so $\det\Phi >0$ and the points $y_0$ and $y_1$  come in
the opposite order.

\vspace{5pt}
The next several sections are therefore devoted to the proof that if
$g_F$ has negative scalar curvature, then the only way in which $g_F$
can extend smoothly to the special orbits, is the `standard way'
described above.

\subsection{Fermi coordinates}
\label{fermi} 
Near any point of $S$ we can introduce {\em Fermi coordinates}. For 
$p$ near $S$, we write  $r(p)$ for the distance from $p$ to $S$ and 
use an angular coordinate $\theta$ in the normal bundle of $S$ in 
$X$. The metric takes the form
$$
dr^2 +r^2d\theta^2 + h_1 + rh_2 + r^2h_3
$$
where $h_1$ is the `first fundamental form' (restriction of the metric 
to $S$), $h_2$ is the second fundamental form, and $h_3$ is a
form on $TX$ bilinear in $rd\theta$ and $TS$. By Gauss's lemma, 
$h_3$ does not contain terms in $dr$, though we shall not use this.

Since $S$ is a fixed-point set of an isometry, the second 
fundamental form $h_2$ is actually zero, but again, we shall not need 
this precision.

In our case, we can choose another angular variable $\phi\in T^2$ as one of 
the coordinates in $S$; complete the set with $y$. Then to leading 
order in $r$, we have, near some given point of $S$,
\begin{equation} \label{1.29.5}
g= dr^2 + dy^2 + r^2d\theta^2 + a^2d\phi^2
\end{equation}
where $a>0$ at $r=0$.  We shall now compare this with the metric $g_F$.

From \S\ref{sec:loc}, if $g_F$ has negative scalar curvature,
\begin{equation}\label{10.29.5}
g_F = \frac{\rho(f_\rho^2+f_\eta^2)-ff_\rho}{\rho f^2}(d\rho^2+d\eta^2)
+
\frac{\rho}{f^2}\frac{1}{\rho(f_\rho^2+f_\eta^2)-ff_\rho}
(d\psi_1,d\psi_2) P^tP\begin{pmatrix} d\psi_1\\ d\psi_2\end{pmatrix}
\end{equation}
where 
$$
P =\begin{pmatrix} \rho f_\eta -\eta f_\rho & f_\rho\cr
f- \rho f_\rho -\eta f_\eta & f_\eta\cr
\end{pmatrix}.
$$
We shall equate the fibre part of this to the `angular part' of
\eqref{1.29.5} and work to leading order in $r$. We shall not divide 
by any quantity tending to zero with $r$ at this stage, so we can 
afford to ignore the higher-order corrections to the metric. Since 
$$
\det P = \rho(f_\rho^2+f_\eta^2)-ff_\rho
$$
we can write this equation as
\begin{equation}\label{2.29.5}
\frac{\rho}{f^2}\frac{1}{\det P}P^tP =
\begin{pmatrix} a^2 & 0\cr 0 & r^2\cr\end{pmatrix} +\cdots.
\end{equation}
It follows that
\begin{equation}
R:= \frac{1}{f}\frac{\sqrt{\rho}}{\sqrt{\det P}}P
\begin{pmatrix} a^{-1} & 0\cr 0 & r^{-1}\cr\end{pmatrix} +\cdots
\end{equation}
is an orthogonal matrix. In fact $\det R = 1$ so we may write
$$
R = \begin{pmatrix} c & s\cr -s & c\cr\end{pmatrix},\;\;c^2+s^2=1.
$$
Hence we find 
\begin{equation}\label{3.29.5}
\frac{\sqrt{\rho}}{f\sqrt{\det P}}
\begin{pmatrix} 
\rho f_\eta -\eta f_\rho & f_\rho\cr
f- \rho f_\rho -\eta f_\eta & f_\eta\cr
\end{pmatrix}
= 
\begin{pmatrix} c a & sr \cr -sa & cr\cr
\end{pmatrix}+\cdots.
\end{equation}
This equation contains most of the information needed to prove 
Theorem~\ref{thmd}. To get our hands on it we proceed as follows: eliminate 
$f_\rho$ and $f_\eta$ between the $(1,1)$, $(1,2)$ and $(2,2)$ 
components of this equation, to obtain
$$
(r\rho-a)c = r\eta s
$$
and hence
\begin{equation*}
c=\frac{r\eta}{\sqrt{(r\rho-a)^2 +r^2\eta^2}}+\cdots,\quad
s=\frac{r\rho-a}{\sqrt{(r\rho-a)^2 +r^2\eta^2}}+\cdots .
\end{equation*}
Next, eliminate $f_\rho$ and $f_\eta$ between the $(1,2)$, $(2,2)$ 
and $(2,1)$ components of this equation, to obtain
$$
\frac{\sqrt{\rho}}{\sqrt{\det P}} -r\rho s - r\eta c = -sa
$$
so that
\begin{equation}\label{5.29.5}
\frac{\sqrt{\rho}}{\sqrt{\det P}}
=\sqrt{(r\rho-a)^2 +r^2\eta^2}+\cdots.
\end{equation}
Finally, we note by taking the determinant of \eqref{3.29.5},
\begin{equation}\label{6.29.5}
\frac{\rho}{f^2} = ar+\cdots.
\end{equation}

\subsection{Proof that the special orbits are at $\rho=0$}

Select a point $p$ on $S$ and suppose that the orbit over $p$ 
corresponds to the point $(\rho_0,\eta_0)$ in the closure of the upper 
half-plane. Suppose for a contradiction that we have $\rho_0>0$.  

It follows from \eqref{5.29.5} that $\sqrt{\rho/\det P}$ is 
bounded away from zero and infinity in the limit.  Hence, from the 
second column of \eqref{3.29.5} it follows that $\del_\rho\log f$ 
and $\del_\eta\log f$ remain bounded on approach to $(\rho_0,\eta_0)$.
Hence $\log f$ is bounded there and so $f$ is bounded away from $0$ 
and $\infty$ there. This now contradicts \eqref{6.29.5} which 
requires $f\to\infty$ on approach to $(\rho_0,\eta_0)$.

\subsection{Boundary behaviour of $f$}

We now know that $\rho\to 0$ as $r\to 0$. Repeating the arguments of 
the previous subsection (but now with $\rho_0=0$), we see that 
$\del_\rho\log f$ and $\del_\eta\log f$ both tend to zero as $\rho\to 
0$. In fact, we obtain
\begin{equation*}
f(\rho,\eta) = A + \frac{1}{2}\rho^2f_1(\eta)+\cdots
\end{equation*}
where $A$ is constant and $f_1$ is some function of
$\eta$. From~\eqref{5.29.5}, we see also that, to leading order in $\rho$,
\begin{equation*}
\frac{\det P}{\rho} = \frac{1}{a^2},
\end{equation*}
and from \eqref{6.29.5}
\begin{equation*}
a\,dr = d\rho/A^2.
\end{equation*}
Hence by comparing coefficients of $dr^2$ in \eqref{1.29.5} and 
\eqref{10.29.5},
$$
1 = \frac{1}{A^2 a^2}(A^2 a)^2 = A^2.
$$
Hence $A=1$ since we are assuming (as we can) that $f>0$ near $\rho=0$.

\subsubsection{Remark} If the isotropy is assumed to be $G(m,n)$ 
then we would have found $f(\rho,\eta) = m\eta -n + 
\frac{1}{2}\rho^2f_1(\eta)$ instead.  That is, we have recovered the
`standard  situation' described at the beginning of this section.

\subsection{Corner behaviour of $f$}

We have now established that $f_0$ has the `standard behaviour' 
\eqref{2.18.3} away from the corners $y_0$ and $y_1$. It remains to 
show that $f_0$ is continuous at these points.  To do this, we repeat 
the discussion above, but with the metric expanded about the point 
$S_0\cap S_1$. It is now natural to
introduce parameters $r_0$ and $r_1$, the 
distance functions from $S_0$ and $S_1$ respectively. Since
$S_0$ and $S_1$ are totally geodesic, the restriction of $r_1$ to 
$S_0$ is also the distance function from the point $S_0\cap S_1$ and 
similarly for $r_0$. Therefore the metric takes the form
\begin{equation*}g =  dr_0^2 + dr_1^2 + r_0^2d\vartheta_0^2 + 
r_1^2d\vartheta_1^2
+\cdots
\end{equation*}
where $\cdots$ denotes terms of higher order.  If we carry through the 
calculations of \S\ref{fermi} with this metric, then we obtain the formulae
\begin{align} 
\frac{\rho}{f^2} & = r_0r_1+\cdots \\
(\log f)_\rho & =
\frac{r_1^2\rho  - r_0r_1}{r_1^2\eta^2 + (r_0-r_1\rho)^2}+\cdots
\label{11.29.5} \\
(\log f)_\eta & =
\frac{r_1^2\eta}{r_1^2\eta^2 + (r_0-r_1\rho)^2}+\cdots \label{12.29.5}
\end{align}
Let the corner correspond to $\eta=\eta_0$. We shall investigate the 
behaviour of $\log f$ in polar coordinates centred at $(0,\eta_0)$, by 
introducing
$$
\eta -\eta_0 = R\cos\Theta,\quad \rho = R\sin\Theta.
$$
Then 
$$
\del_R\log f  = (\cos\Theta)\, \del_\eta\log f +(\sin\Theta)\,\del_\rho\log f
$$
and this is uniformly bounded for any fixed $\Theta$ by \eqref{11.29.5} 
and \eqref{12.29.5}. Hence $f(R,\Theta)$ has a limit as $R\to 0$, for 
each fixed $\Theta$. To see that these limits are independent of 
$\Theta$, we note that
$$
\del_\Theta\log f = R(-\sin\Theta\, \del_\eta + \cos\Theta \,\del_\rho)\log f
$$
and the right hand side is $O(R)$, uniformly in $\Theta$. Hence
$$
|f(R,0) -f(R,\pi)| = O(R)
$$
and this shows that $f$ is continuous at the corner, by taking $R\to
0$. In view of the discussion in \S\ref{outline}, the proof of
Theorem~\ref{thmd} is now complete.

\section*{Appendix}
In this appendix we explain the geometric origin of the basic solutions
that we have used in this paper. Much of the discussion goes through for
hyperbolic space of arbitrary dimension, but we shall confine ourselves to
$\cH^2$. As in \cite{CaPe:emt}, we shall fix a $2$-dimensional symplectic
vector space $\W$ and consider the symmetric square $S^2\W$, equipped with
the quadratic form $\det$ as a 3-dimensional Minkowksi space. If $g\in
SL(\W)$, $x\in S^2\W$, then $g\cdot x = gxg^t$ is evidently an isometric
action, corresponding to the double cover of the identity component of
$SO(1,2)$ by $SL_2(\R)$.

The hyperbolic plane $\cH^2$ appears as (one sheet of) the hyperboloid
$\det x =1$ or essentially equivalently as the open subset $S_+/\R^+$ of
the projective space $P(S^2\W)$, where
\begin{equation*}
S_+ = \{x\in S^2\W: \det x >0\}.
\end{equation*}
It is easy to check that there is a natural 
correspondence
\begin{equation}\label{eq:efns}
\{f\in \Cinf(S_+): \Quabla f=0,\ E\cdot f = \alpha f\}
= \{g\in \Cinf(\cH^2): \Delta g = \alpha(\alpha + 1)g\},
\end{equation}
where $E$ is the Euler homogeneity operator on $S^2\W$ and $\Quabla$
is the wave operator of $S^2\W$.  If $N$ is any nonzero null vector
($\det N=0$), then $x\mapsto f(N\cdot x)$ is a solution of the wave equation,
for any function $f$.  The `basic solutions' of the equation 
$\Delta F = (3/4)F$ arise in precisely this way, with the simplest
possible function homogeneous of degree $1/2$, namely $(N\cdot x)^{1/2}$.

If we represent $N = n\otimes n$ where $n=(1,y)$ and use the
parameterization of the hyperboloid by half-space coordinates
\begin{equation*}
x= \frac{1}{\rho}\begin{bmatrix} 1 & \eta \cr
\eta & \rho^2+\eta^2\cr\end{bmatrix}
\end{equation*}
then 
\begin{equation*}
\sqrt{N\cdot x}=\sqrt{n x^{-1} n^t}
=\frac{\sqrt{\rho^2 + (\eta -y)^2}}{\sqrt{\rho}}
\end{equation*}
which is the `basic solution' $F(\rho,\eta;y)$ of \eqref{27.1.22}.

Just as eigenfunctions of the laplacian in $\cH^2$ correspond to
homogeneous solutions of the wave equation in $S_+$, so eigenfunctions of
the Dirac operator in $\cH^2$ correspond to homogeneous solutions of the
Dirac equation in $S_+$. Here we can identify the product bundle
$S_+\times\W$ as the spin-bundle of $S_+$; its restriction to the
hyperboloid is naturally isomorphic to the spin-bundle of $\cH^2$ (though
the induced metric at $x$ is given by $x^{-1}$).  It is easy to check that
if $N = n\otimes n$ as before, then the $\W$-valued function
$n/\sqrt{N\cdot x}$ is a solution of the minkowskian Dirac equation. This
descends to the `basic solution' of the Joyce equation $\phi(\rho,\eta;y)$
\eqref{27.1.8} after replacing $n$ by the dual vector $n^*$, $\ip{n,u} =
n^*(u)$. (Explicitly, if $n = (1,y)$ then $n^* = (-y,1)$.)

Over the conformal boundary $\RP{1}$ of $\cH^2$ there is a family of
homogeneous line bundles. In terms of the parameter $n\in \W$, we define
the bundle $\cO(1)$ as corresponding to functions homogeneous of degree $1$,
\begin{equation*}
f(\lambda n) = \lambda f(n), \quad \lambda\neq 0,
\end{equation*}
and $|\cO(1)|$ to correspond to functions homogeneous in the sense
\begin{equation*}
f(\lambda n) = |\lambda| f(n), \quad \lambda\neq 0.
\end{equation*}
Notice that $\cO(1)=|\cO(1)|\tens L$ where $L$ is the M\"obius bundle with
$L^2=\cO$, and that $|\cO(1)|$ is topologically trivial, whereas $\cO(1)$ is
not.  These considerations are important in trying to find the correct
interpretation of the formulae for the smeared solutions in
\S\ref{sec:smear}.  If $M\to \RP{1}$ is a given line bundle, then by a
distributional section of $M$ we mean a continuous linear functional on the
space $\Cinf(\RP{1},M^*\otimes \Omega)$. Here $\Omega$ is the bundle of
densities, which can be identified invariantly with $\cO(-2)$.  In so far
as $\sqrt{N\cdot x}$ is a section of $\cO(1)\tens L$ for each fixed $x$, it
follows that if $u$ is a distributional section of $\cO(-3)\otimes L$ then
\begin{equation*}
\int u(n)\sqrt{N\cdot x}
\end{equation*}
is a well-defined function of $x$. This is the $SL_2(\R)$-invariant
interpretation of formula \eqref{4.1.2}. The formula \eqref{3.1.2} can
be interpreted similarly.

\nocite{BGMR:3s7}
\nocite{GiHa:gmi}
\nocite{Bel:ncr}
\nocite{Rol:reh}
\nocite{LeBr:cqk}
%
% Bibliography format
%
\newcommand{\bauth}[1]{\mbox{#1}} \newcommand{\bart}[1]{\textit{#1}}
\newcommand{\bjourn}[4]{#1\ifx{}{#2}\else{ \textbf{#2}}\fi{ (#4)}}
\newcommand{\bbook}[1]{\textsl{#1}}
\newcommand{\bseries}[2]{#1\ifx{}{#2}\else{ \textbf{#2}}\fi}
\newcommand{\bpp}[1]{#1} \newcommand{\bdate}[1]{ (#1)} \def\band/{and}
\newif\ifbibtex%\bibtextrue% Comment out \bibtextrue to remove BiBTeX.
\ifbibtex
\bibliographystyle{genbib}
\bibliography{papers}
\else

\fi

\end{document}